\documentclass[a4paper,10pt]{article}
\usepackage{fontenc}

\usepackage{amsmath}
\usepackage{amssymb}
\usepackage{amsthm}
\usepackage{geometry}
\usepackage{graphicx}
\usepackage{times}
\usepackage{xcolor}
\usepackage{etoolbox}

\date{}

\newtoggle{arxiv}
\toggletrue{arxiv}
\iftoggle{arxiv} 
{
  \usepackage{hyperref}
}
{
  \usepackage[colorlinks]{hyperref}
  \usepackage{microtype}
  \DeclareOldFontCommand{\rm}{\normalfont\rmfamily}{\mathrm}
}

\theoremstyle{plain}
\newtheorem{Theorem}{Theorem}[section]
\newtheorem{Proposition}[Theorem]{Proposition}

\newtheorem{Corollary}[Theorem]{Corollary}
\newtheorem{Lemma}[Theorem]{Lemma}

\newtheorem{Observation}[Theorem]{Observation}

\theoremstyle{definition}

\theoremstyle{remark}

\newcommand{\gauss}[3]{\genfrac{[}{]}{0pt}{}{#1}{#2}_{#3}}
\newcommand{\PG}[2]{\operatorname{PG}(#1,#2)}
\newcommand{\F}[2]{\mathbb{F}_{#2}^{#1}}
\newcommand{\GF}{\mathbb{F}}
\newcommand{\G}[3]{\gauss{\mathbb{F}_{#3}^{#1}}{#2}{}}

\begin{document}
\title{Tables of subspace codes}
\author{Daniel~Heinlein, Michael~Kiermaier, Sascha~Kurz, and Alfred~Wassermann\thanks{
All authors are or were with the Department of Mathematics, Physics, and Computer Science, University of Bayreuth, Bayreuth, GERMANY. email: firstname.lastname@uni-bayreuth.de\newline
The work was supported by the ICT COST Action IC1104 and grants KU 2430/3-1, WA 1666/9-1 -- ``Integer Linear Programming Models for Subspace Codes and Finite Geometry'' -- from the German Research Foundation.}}
\maketitle
\begin{abstract}
One of the main problems of subspace coding asks for the maximum possible cardinality of a subspace code
with minimum distance at least $d$ over $\F{n}{q}$, where the dimensions of the codewords, which are vector spaces, 
are contained in $K\subseteq\{0,1,\dots,n\}$. In the special case of $K=\{k\}$ one speaks of 
constant dimension codes. Since this (still) emerging field is very prosperous on the one hand side and 
there are a lot of connections to classical objects from Galois geometry it is a bit difficult 
to keep or to obtain an overview about the current state of knowledge. To this end we have 
implemented an on-line database of the (at least to us) known results at \url{subspacecodes.uni-bayreuth.de}.
The aim of this technical report is to provide a user guide how this technical 
tool can be used in research projects and to describe the so far implemented theoretic and algorithmic 
knowledge.   

\medskip

\end{abstract}

\section{Introduction}

\noindent
The seminal paper by K\"otter and Kschischang \cite{koetter2008coding} 
started the interest in subspace codes which are sets of 
subspaces of the $\GF_q$-vector space $\GF_q^n$.
Two widely used distance measures for subspace codes (motivated 
by an information-theoretic analysis of the 
K\"otter-Kschischang-Silva model, see e.g.\ \cite{silva2008rank}) are the 
\emph{subspace distance} 
\[
  d_S(U,W):=\dim(U+W)-\dim(U\cap W)=2\cdot\dim(U+W)-\dim(U)-\dim(W)
\] 
and the
\emph{injection distance} 
\[
  d_I(U,W):=\max\left\{\dim(U),\dim(W)\right\}-\dim(U\cap W), 
\]
where $U$ and $W$ are subspaces of $\GF_q^n$.
The two metrics are equivalent, i.e., it is known that $d_I(U,W) \le d_S(U,W) \le 2d_I(U,W)$.
Here, we restrict ourselves to the subspace distance.

The set of all $k$-dimensional subspaces of an $\GF_q$-vector space $V$ will be denoted by $\gauss{V}{k}{q}$.
For $n = \dim(V)$, its cardinality is given by the Gaussian binomial coefficient
\[
\gauss{n}{k}{q} =
\begin{cases}
  \frac{(q^n-1)(q^{n-1}-1)\cdots(q^{n-k+1}-1)}{(q^k-1)(q^{k-1}-1)\cdots(q-1)} & \text{if }0\leq k\leq n\text{;}\\
  0 & \text{otherwise.}
\end{cases}
\]

A set $\mathcal{C}$ of subspaces of $V$ is called a \emph{subspace code}.
The \emph{minimum distance} of $\mathcal{C}$ is given by $d = \min\{d_S(U,W) \mid U,W\in\mathcal{C}, U \neq W\}$.
If the dimensions of the \emph{codewords}, i.e., the elements of $\mathcal{C}$ are contained in some set 
$K \subseteq \{1,\ldots,n\}$, $\mathcal{C}$ is called an $(n,\#\mathcal{C},d;K)_q$ subspace code.
In the unrestricted case $K = \{0,\ldots,n\}$, also called mixed dimension case, we use the notation $(n,\#\mathcal{C},d)_q$ subspace code.
In the other extreme case $K = \{k\}$, we use the notation $(n,\#\mathcal{C},d;k)_q$ and call $\mathcal{C}$ a \emph{constant dimension code}.

For fixed ambient parameters $q$, $n$, $K$ and $d$, a \emph{main problem of subspace coding} asks for the determination of the maximum 
possible size $A_q(n,d;K) := M$ of an $(n,M,\ge d;K)_q$ subspace code and~--~as a refinement -- the classification of all corresponding 
optimal codes up to isomorphism.
Again, the simplified notations $A_q(n,d)$ and $A_q(n,d;k)$ are used for the unrestricted case $K = \{0,\ldots,n\}$ and the constant 
dimension case $K = \{k\}$, respectively. Note that in the latter case $d_S(U,W)=2\cdot d_I(U,W)\in 2\cdot\mathbb{N}$ is an even number.

In general, the exact determination of $A_q(n,d;K)$ is a hard problem, both on the theoretic and the algorithmic side.
Therefore, lower and upper bounds on $A_q(n,d;K)$ have been intensively studied in the last years, see e.g.\ \cite{etzion2013problems,heinlein2018binary,
MR3543542}. Since the underlying discrete structures arose under different names in different fields of discrete mathematics, it is even more difficult 
to get an overview of the state of the art. 
For example, geometers are interested in so-called partial $(k-1)$-spreads of $\PG{n-1}{q}$. Following 
the track of partial spreads, one can end up with orthogonal arrays or $(s,r,\mu)$-nets.
Furthermore, $q$-analogs of Steiner systems provide optimal constant dimension codes.
For some sets of parameters constant dimension codes are in one-to-one correspondence with so-called vector space partitions. 

The aim of this report is to describe the underlying theoretical base of an 
on-line database, found at

\noindent
\begin{center}
\href{http://subspacecodes.uni-bayreuth.de}{http://subspacecodes.uni-bayreuth.de}
\end{center}

\noindent
and maintained by the authors that
tries to collect up-to-date information 
on the best lower and upper bounds for subspace codes. 
Whenever the exact value $A_q(n,d;K)$ could be determined, 
we ask for a complete classification of all optimal codes up to isomorphism. 
Occasionally we list classifications for non-maximum codes, too.
Since the overall task is rather comprehensive, 
we start by focusing on the special cases of constant dimension codes, $A_q(n,d;k)$, and (unrestricted) 
subspace codes, $A_q(n,d)$, using the subspace distance as metric.
For a more comprehensive survey on network coding we refer the interested reader e.g.\ to \cite{bassoli2013network}.
For algorithmic aspects we refer the interested reader e.g.\ to \cite{kohnert2008construction}.

The remaining part of this report is structured as follows. In Section~\ref{sec_structure_website} 
we outline the structure of the website and how to access the data. Constant dimension codes (CDC) are treated in Section~\ref{sec_bounds_cdcs}, 
where the currently implemented lower bounds, constructions, and upper bounds are 
described in Subsection~\ref{subsec_implemented_constructions} and Subsection~\ref{subsec_cdc_upper}, 
respectively. Mixed dimension codes (MDC) are treated in Section~\ref{sec_mdc}, where the implemented lower bounds,
constructions, and upper bounds are described in Subsection~\ref{subsec_mdc_lower} and Subsection~\ref{subsec_mdc_upper},  
respectively. Finally we draw a conclusion in Section~\ref{sec_conclusion} 
and list some explicit tables on upper and lower bounds in an appendix.

\section{Structure of the website}
\label{sec_structure_website}

\begin{figure}[ht]
\iftoggle{arxiv}
{\includegraphics[width=\textwidth]{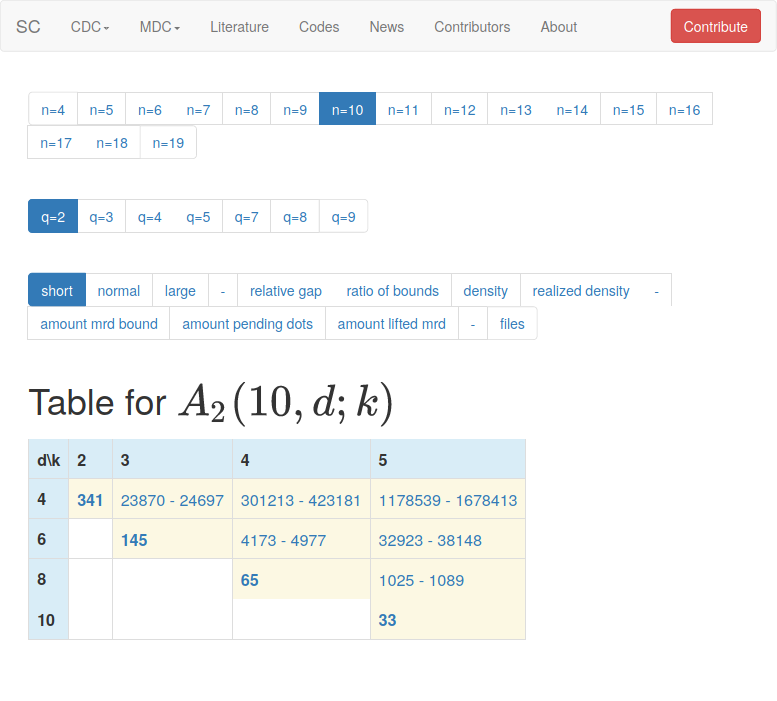}}
{\includegraphics[width=\textwidth]{images/cdc_2019}}
\caption{Tables of constant dimension codes}\label{pic:cdc}
\end{figure}

\noindent
On the website the two special cases $A_q(n,d;k)$ and $A_q(n,d)$ can be accessed via the menu items
\texttt{CDC} (constant dimension code) and \texttt{MDC} (mixed dimension code), see Figure \ref{pic:cdc}. Selecting the item \texttt{Table} 
yields the rough data that we will outline in this section. Selecting the item \texttt{Constraints} yields information 
about the so far implemented lower and upper bounds.

\subsection{Constant dimension codes -- \texttt{CDC}}
\label{subsec_table_cdc}

\noindent
For a constant dimension code the dimension $n$ of the ambient space (first \textit{selection row}) and the 
field size $q$ (second \textit{selection row}) can be chosen. The current 
limits are $2\le q\le 9$ and $4\le n\le 19$ (resp. in the \emph{large} view $1 \le n \le 19$). For each chosen pair of those parameters a table with the 
information on lower and upper bounds on constant dimension codes over $\F{n}{q}$ is displayed.

The rows of those tables are labeled by the minimum distance $d=d_S(\star)$ and the columns are labeled by the 
dimension $k$ of the codewords. In the third \textit{selection row} several \textit{views} can be picked. The first three 
options, \texttt{short}, \texttt{normal}, and \texttt{large}, specify the subset of possible values for the parameters $d$ 
and $k$.  In the most extensive view \texttt{large}, $k$ can take all integers between $0$ and $n$. For $d$ the integers 
between $1$ and $n$ are considered. As
\begin{itemize}
  \item $A_q(n,d;0)=1$ for all $1\le d\le n$;
  \item $A_q(n,d;k)=A_q(n,d;n-k)$;
  \item $A_q(n,2d'+1;k)=A_q(n,2d'+2;k)$ for all $d'\in\mathbb{N}$;
\end{itemize} 
one may assume $1\le k\le\left\lfloor n/2\right\rfloor$, $2\le d\le n$, and $d\in 2\mathbb{N}$. These assumptions 
are implemented in the view \texttt{normal}. However, some exact values of $A_q(n,d;k)$ are rather easy to determine
\begin{itemize}
  \item $A_q(n,2;k)=\gauss{n}{k}{q}$, since any two different $k$-dimensional subspaces of $\F{n}{q}$ have a subspace 
        distance of at least $2$;
  \item if $d>2k$, then we can have at most one codeword, i.e., $A_q(n,d;k)=1$.      
\end{itemize} 
Thus, we may assume $2\le k\le\left\lfloor n/2\right\rfloor$, $4\le d\le 2k$, and $d\in 2\mathbb{N}$. These assumptions 
are implemented in the view \texttt{short}. The standard selection is given by $n=4$, $q=2$ and the view \texttt{short}.

Given one of these three views, a table entry may consist of
\begin{itemize}
  \item a range $l$--$u$: An example is given by the parameters $q=2$, $n=7$, $d=4$, $k=3$, where $l=333$ and $u=381$.
        The meaning is that for the corresponding maximum cardinality of a constant dimension code only the lower bound $l$ 
        and the upper bound $u$ is known, i.e., $333\le A_2(7,4;3)\le 381$ in the example.\footnote{The lower bound $333$ has been 
        determined in \cite{heinlein2019subspace}. If the upper bound $381$ is attained, then the corresponding automorphism group 
        can have an order of at most two~\cite{fano_aut}.}
  \item a \textbf{bold} number $m$: An example is given by the parameters $q=2$, $n=10$, $d=8$, $k=4$, where $m=65$.
        The meaning is that the corresponding maximum cardinality of a constant dimension code is exactly determined, i.e., 
        $A_2(10,8;4)=65$ in the example.
  \item a \textbf{bold} number $m$ with an asterisk and a number $l$ in brackets: An example is given by the parameters $q=2$, 
        $n=6$, $d=4$, $k=3$, where $m=77$ and $l=5$. The meaning is that the corresponding maximum cardinality of a constant 
        dimension code is exactly determined and all optimal codes have been classified up to isomorphism, i.e., 
        $A_2(6,4;3)=77$ and there are exactly $5$ isomorphism types in the example, see \cite{hkk77}. Another example is given for the 
        parameters $q=2$, $n=6$, $d=4$, and $k=2$, where there are exactly $131,044$ isomorphism types of constant dimension 
        codes attaining cardinality $A_2(6,4;2)=21$, see \cite{MR2475427}.
  \item a \textbf{bold} number $m$ with a lower bound $\ge l$ in brackets: An example is given by the parameters $q=2$, 
        $n=13$, $d=4$, $k=3$, where $m=1597245$ and $l=512$. The meaning is that the corresponding maximum cardinality of a constant 
        dimension code is exactly determined and there are at least $l$ isomorphism classes of optimal codes.
\end{itemize}
Each nontrivial table entry is clickable and then yields further information on several lower and upper bounds, see 
Subsection~\ref{subsec_implemented_constructions} and Subsection~\ref{subsec_cdc_upper} for the details. 

In some cases, e.g., for the parameters $q=2$, $n=6$, $d=4$, and $k=3$, the corresponding codes are also available 
for download using the button called ``file''. The format of these codes is mostly GAP\footnote{\url{http://www.gap-system.org}} or MAGMA\footnote{\url{http://magma.maths.usyd.edu.au}}.

Besides the views \texttt{short}, \texttt{normal}, and \texttt{large} for the selection of ranges for the parameters 
$d$ and $k$, there are some additional views. The views \texttt{relative gap} and \texttt{ratio of bounds} condense the current  
lack of knowledge on the exact value of $A_q(n,d;k)$ to a single number. For the view \texttt{relative gap} this number is 
given by the formula
\[
  \frac{\text{upper bound}\,-\,\text{lower bound}}{\text{lower bound}},
\]
i.e., we obtain a non-negative real number. While principally any number in $\mathbb{R}_{\ge 0}$ can be obtained, 
the largest relative gap in our database is currently given by about $0.619$ for the parameters $q=2,n=19,d=4,k=9$. A gap of $0.0$ 
corresponds to the determination 
of the exact value $A_q(n,d;k)$. The mentioned formula is also displayed on the webpage, when you move your mouse 
over the word \texttt{relative gap}. For the view \texttt{ratio of bounds} the corresponding number is given by the formula
\[
  \frac{\text{lower bound}}{\text{upper bound}},
\] 
which may take any real number in $(0,1]$. The smallest ratio of bounds in our database is given by about $0.618$ for the same parameters as above. 
Clearly, the largest relative gap yields the smallest ratio of bounds and vice versa as the function $x \mapsto \frac{1}{x}-1$ is strictly decreasing in $(0,1]$. A ratio of bounds of $1.0$ corresponds to the determination of the exact value $A_q(n,d;k)$. The mouse-over effect 
is also implemented in that case.

The views \texttt{density} and \texttt{realized density} compare the Anticode bound, see (Theorem~\ref{theo:anticode}), to the best known upper 
bound and best known lower bound, respectively, i.e.,

\[
\frac{\text{best known upper bound}}{\text{Anticode bound}},
\qquad
\text{and}
\qquad
\frac{\text{best known lower bound}}{\text{Anticode bound}}.
\]

Hence, they are a measure how dense it is possible to fill the Grassmannian with codewords. Note that in the case of Steiner Systems, both bounds, the \texttt{density} and the \texttt{realized density}, are one since the size of a Steiner System is exactly the size of the Anticode bound.

Another type of view arose from some of the various constructions described in Subsection~\ref{subsec_implemented_constructions}. 
They are labeled as \texttt{amount pending dots} and \texttt{amount lifted mrd} and condense the \textit{strength} of a 
certain construction to a single number in $\mathbb{R}_{\ge 1}$. This number is always given as the quotient between the currently best known 
lower bound and the value obtained by the respective construction. Here, a value of one means that the currently best known code 
can be obtained by the respective construction. A value larger than $1$ measures how much better a more tailored construction 
is for this specific set of parameters compared to the respective general construction method. We remark that 
\texttt{amount pending dots} is still experimental and in some cases there may still be better codes obtained from the underlying 
very general construction technique, which has quite some degrees of freedom. With respect to upper bounds the additional 
view \texttt{amount mrd bound} is introduced. Here the displayed single number is given by the currently best known lower bound
divided by the so-called MRD bound, see Subsection~\ref{subsec_mrd_bound}.  

The view \texttt{files} is like the view \texttt{short} but the background gets a green color if there is a downloadable file for these parameters.

\subsubsection{Toplist}

\begin{figure}[ht]
\iftoggle{arxiv}
{\includegraphics[width=\textwidth]{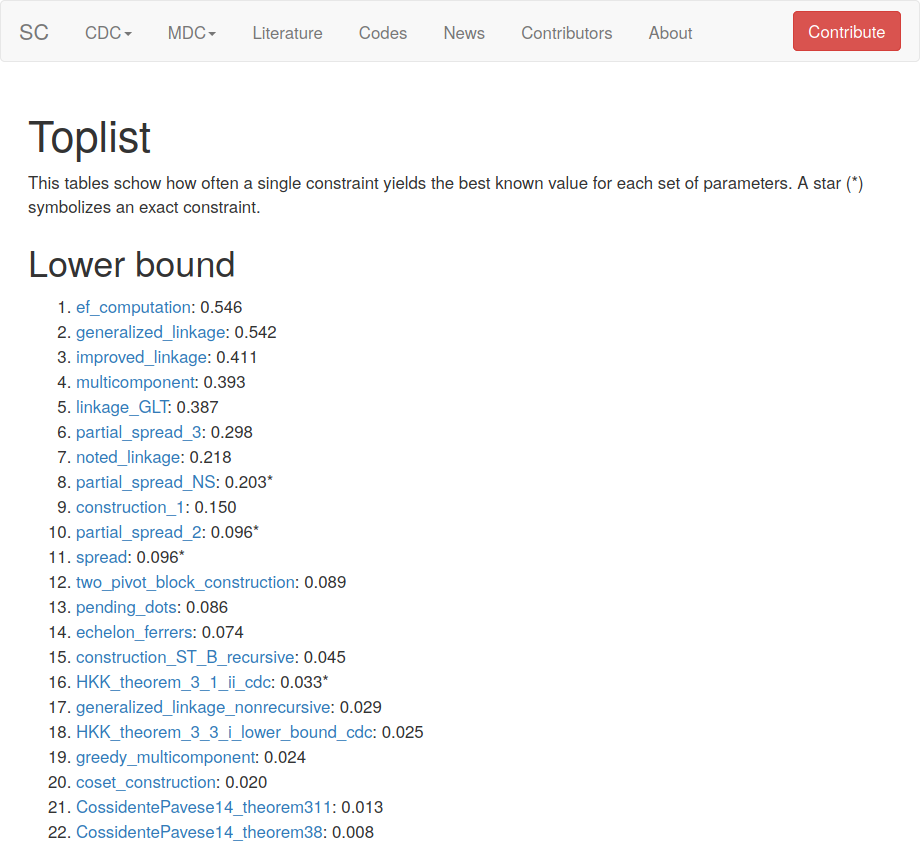}}
{\includegraphics[width=\textwidth]{images/cdc_top_2019}}
\caption{Toplist of constant dimension codes}\label{pic:cdc_top}
\end{figure}

These statistics, see Figure~\ref{pic:cdc_top}, show how often a single constraint yields the best known bound for the parameters $2 \le q \le 9$, $4 \le n \le 19$, 
$2 \le k \le \lfloor n/2 \rfloor$, and $4 \le d \le 2k$, where $d$ is even. For each set of parameters in which a single constraint yields 
the best known value, it scores a point. This score is then divided by the size of the set of parameters, i.e., all constant dimension code 
parameters in the database. Constraints are grouped into two categories: lower and upper bounds and then ordered by their normalized score. 
The special constraints that yield the exact code sizes appear in both categories and are denoted with an asterisk~(*).

Currently the lower bound with the highest score is \texttt{ef\_computation},  
which yields the best known lower bound in 54.6\% of the constant dimension code parameters of the database.

The upper bound \texttt{improved\_johnson} has currently the highest score for upper bounds with 79.0\%.

\subsubsection{Views for single CDCs}
\begin{figure}[ht]
\iftoggle{arxiv}
{\includegraphics[width=\textwidth]{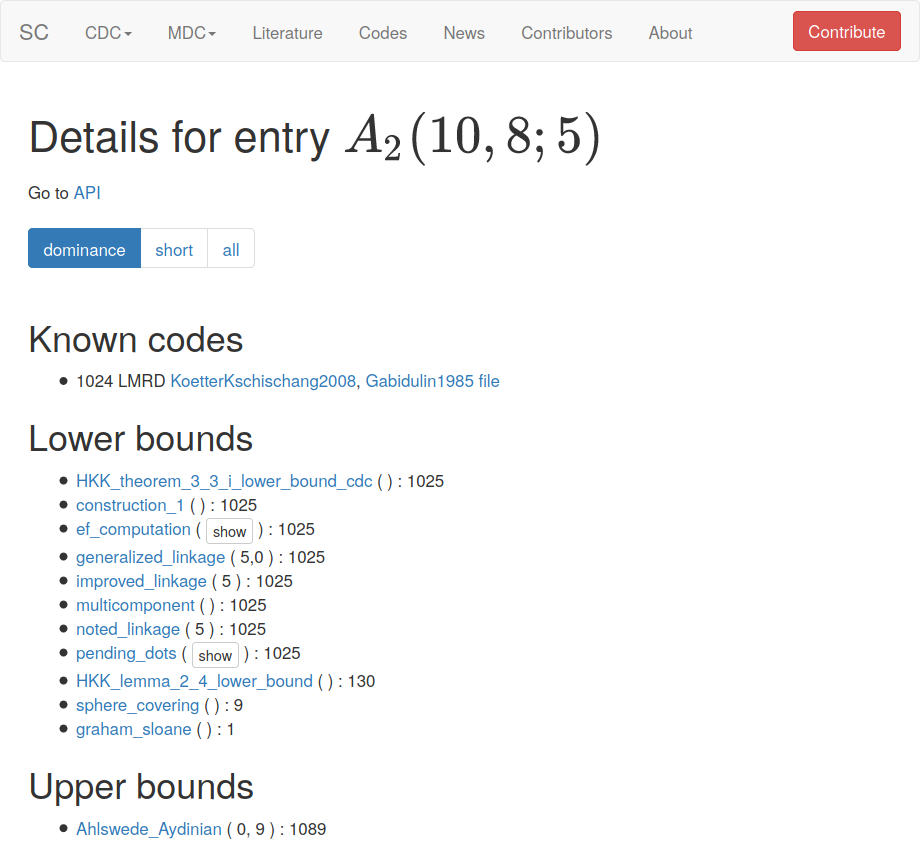}}
{\includegraphics[width=\textwidth]{images/cdc_te_2019}}
\caption{Views of constant dimension code entries}\label{pic:cdc_te}
\end{figure}

Each constant dimension code entry provides multiple level of details which are also called views, see Figure~\ref{pic:cdc_te}.

The view \texttt{all} shows all implemented constraints and the corresponding upper bounds. Some constraints involve a parameter. 
To this end, the view \texttt{short}, which is the default, displays only the best instances for the same constraint (by optimizing the parameter(s)). 
The third view, \texttt{dominance}, performs like \texttt{short} but incorporates an additional filtering due to known relationships between constraints. Hence, if a constraint is always worse than another constraint that is present, it is omitted.

Currently \texttt{dominance} respects for the parameters of this database the following relations between upper bounds:\\
\texttt{sphere\_packing} $\le$ \texttt{all\_subs}\\
\texttt{anticode} $\le$ \texttt{sphere\_packing}\\
\texttt{anticode} $\le$ \texttt{singleton}\\
\texttt{johnson\_1} $\le$ \texttt{johnson\_2}\\
\texttt{johnson\_1} $\le$ \texttt{anticode}\\
\texttt{johnson\_1} $\le$ \texttt{ilp\_1}\\
\texttt{ilp\_1} $\le$ \texttt{ilp\_2}\\
\texttt{ilp\_4} $\le$ \texttt{ilp\_3}\\
\texttt{johnson\_2} $\le$ \texttt{ilp\_4}\\
\texttt{Ahlswede\_Aydinian} $\le$ \texttt{johnson\_1}\\
\texttt{Ahlswede\_Aydinian} $\le$ \texttt{johnson\_2}\\

And between lower bounds:\\
\texttt{sphere\_covering} $\le$ \texttt{trivial\_1}\\
\texttt{echelon\_ferrers} $\le$ \texttt{lin\_poly}\\
\texttt{ef\_computation} $\le$ \texttt{echelon\_ferrers}\\
\texttt{improved\_linkage} $\le$ \texttt{linkage\_GLT}\\
\texttt{improved\_linkage} $\le$ \texttt{linkage\_ST}\\

\subsection{Mixed dimension codes -- \texttt{MDC}}
\label{subsec_table_mdc}

\begin{figure}[ht]
\iftoggle{arxiv}
{\includegraphics[width=\textwidth]{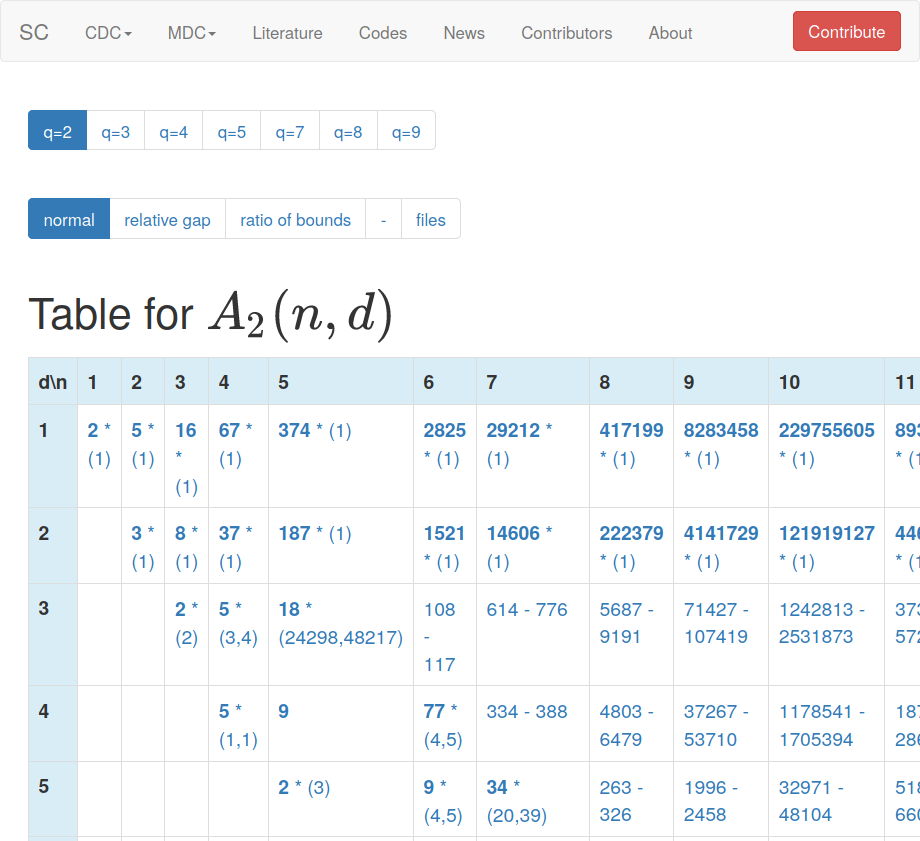}}
{\includegraphics[width=\textwidth]{images/mdc_2019}}
\caption{Tables of (mixed dimension) subspace codes}\label{pic:mdc}
\end{figure}

\noindent
For a subspace code with mixed dimensions the 
field size $q$ (\textit{selection row} number one) can be chosen. The current 
limits are given by  $2\le q\le 9$. For each chosen parameter a table with the 
information on lower and upper bounds on $A_q(n,d)$ over $\F{n}{q}$ ($n \le 19$) is displayed, see Figure \ref{pic:mdc}.

The rows of those tables are labeled by the distance $d=d_S(\star)$ and the columns are label by the 
dimension $n$ of the ambient space $\F{n}{q}$. In the second \textit{selection row} several \textit{views} can be picked. 
The view \texttt{normal}, c.f.~Subsection~\ref{subsec_table_cdc}, already incorporates the restriction to $1 \le d\le n\le 19$.
The views \texttt{relative gap} and \texttt{ratio of bounds} condense the current  
lack of knowledge on the exact value of $A_q(n,d)$ to a single number. For the view \texttt{relative gap} this number is 
given by the formula
\[
  \frac{\text{upper bound}\,-\,\text{lower bound}}{\text{lower bound}},
\]
i.e. we obtain a non-negative real number. While principally any number in $\mathbb{R}_{\ge 0}$ can be obtained, 
the largest relative gap in our database is currently given by about $2.493$ for the parameters $q=2,n=19,d=4$. A relative gap of $0.0$ 
corresponds to the determination 
of the exact value $A_q(n,d)$. The mentioned formula is also displayed on the webpage, when you move your mouse 
over the word \texttt{relative gap}. For the view \texttt{ratio of bounds} the corresponding number is given by the formula
\[
  \frac{\text{lower bound}}{\text{upper bound}},
\] 
which may take any real number in $(0,1]$. The smallest ratio of bounds in our database is given by about $0.29$ for the same parameters as above. Clearly the largest relative gap yields the smallest ratio of bounds and vice versa as the function $x \mapsto \frac{1}{x}-1$ is strictly decreasing in $(0,1]$. A ratio of bounds of $1.0$ corresponds to the determination of the exact value $A_q(n,d)$. The mouse-over effect is also implemented in that case.

The view \texttt{files} is like the view \texttt{normal} but the background gets a green color if there is a downloadable file for these parameters.

\begin{figure}[ht]
\iftoggle{arxiv}
{\includegraphics[width=\textwidth]{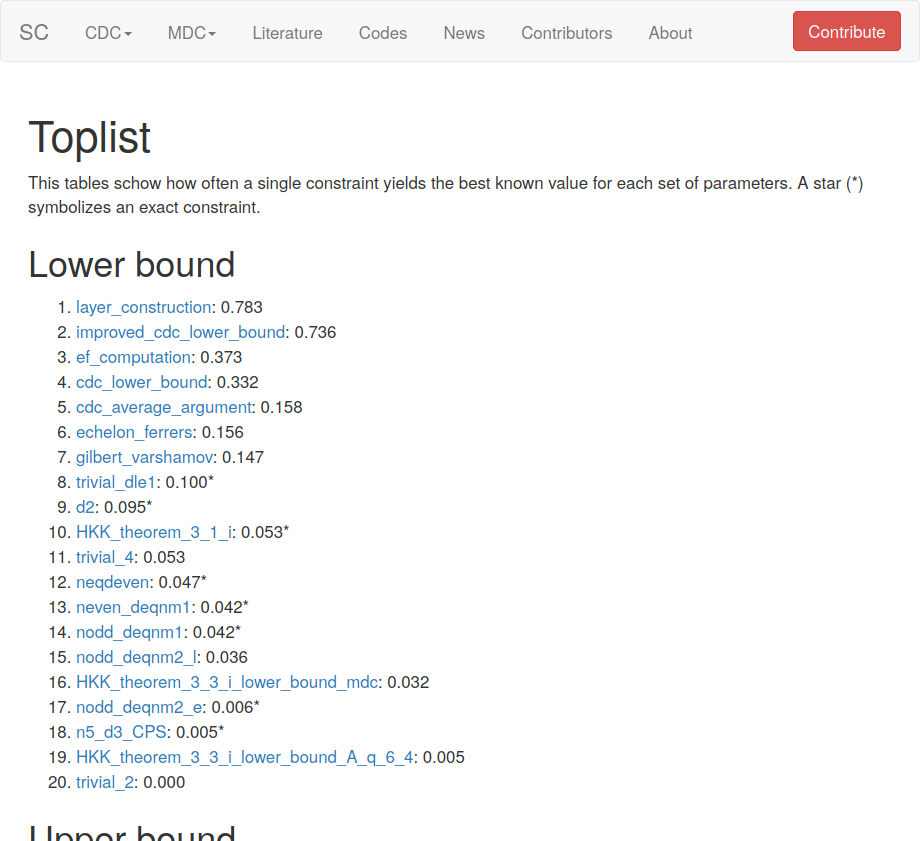}}
{\includegraphics[width=\textwidth]{images/mdc_top_2019}}
\caption{Toplist of (mixed dimension) subspace codes}\label{pic:mdc_top}
\end{figure}

\subsubsection{Toplist}

These statistics show how often a single constraint yields the best known bound for the parameters $2 \le q \le 9$, $4 \le n \le 19$, and $1 \le d \le n$. For each set of parameters in which a single constraint yields the best known value, it scores a point. This score is then divided by the size of the set of parameters, i.e., all mixed dimension code parameters in the database. Constraints are grouped into two categories: lower and upper bounds and then ordered by their normalized score. The special constraints that yield the exact code sizes appear in both categories and are denoted with an asterisk (*).

Currently the lower bound with the highest score is \texttt{layer\_construction} (Theorem~\ref{theo:improved_cdc_lower_bound}) and it yields the best known lower bound in 78.3\% of the mixed dimension code parameters of the database.

The upper bound \texttt{improved\_cdc\_upper\_bound} (Theorem~\ref{theo:improved_cdc_upper_bound}) has currently the highest score for 
upper bounds with 42.3\%, see Figure~\ref{pic:mdc_top}.

\subsubsection{MDC table for arbitrary \texorpdfstring{$q$}{q}}

\begin{figure}[ht]
\iftoggle{arxiv}
{\includegraphics[width=\textwidth]{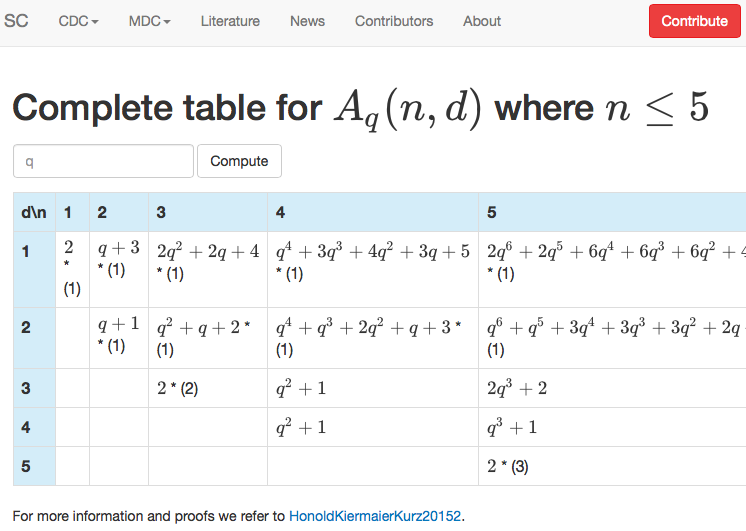}}
{\includegraphics[width=\textwidth]{images/mdc_q}}
\caption{Mixed dimension subspace codes for arbitrary $q$}\label{pic:mdc_q}
\end{figure}

This table has the same layout as the other MDC tables. The benefit is that the optimal sizes for $A_q(n,d)$ is known for all $q$ and 
$n \le 5$, as well as the number of isomorphism types in some of these cases. This table, see Figure~\ref{pic:mdc_q}, show the sizes in 
terms of $q$-polynomials that may even be evaluated for specific $q$ values by entering them in the input box and pressing 
the ``Compute'' button or the ``Enter'' key. It is possible to leave the input field blank to get the $q$-polynomials back or to enter 
non prime power, even negative, numbers.

\subsection{Downloadable files -- \texttt{Codes}}

\begin{figure}[ht]
\iftoggle{arxiv}
{\includegraphics[width=\textwidth]{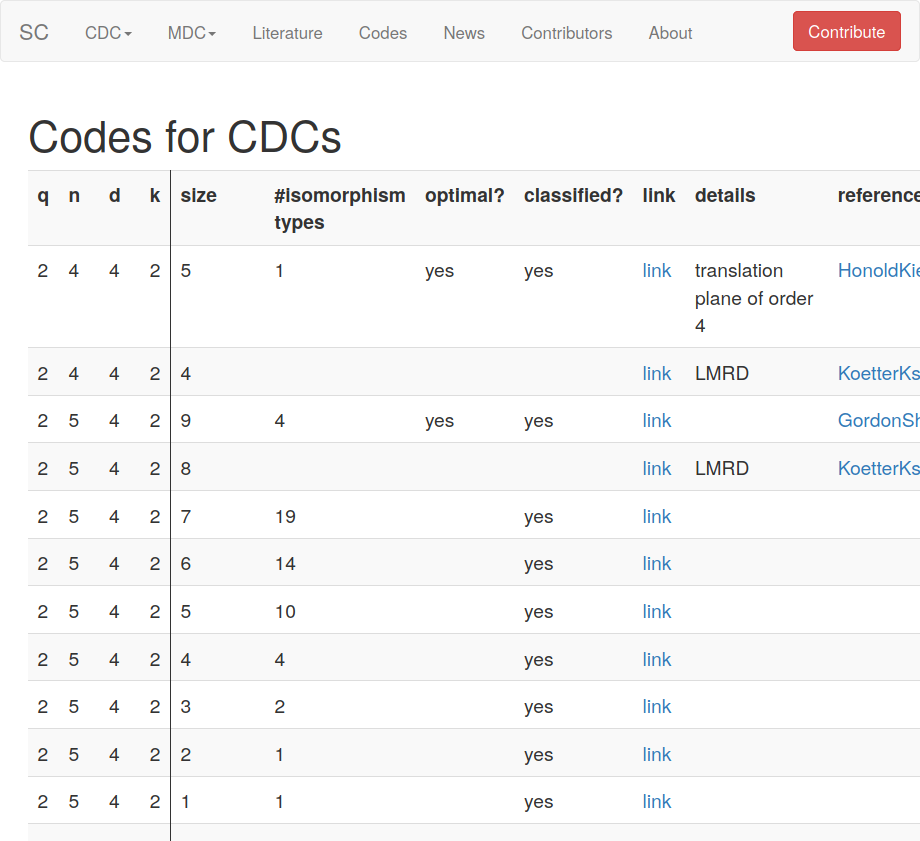}}
{\includegraphics[width=\textwidth]{images/codes_files}}
\caption{A tables of downloadable files}\label{pic:codes_files}
\end{figure}

Under the heading \texttt{Codes} we provide a summarizing table of all codes for which a 
file can be downloaded, see Figure~\ref{pic:codes_files}. The information includes the parameters $q$, $n$, $d$, $k$, the size of 
the code, the number of isomorphism types if the codes are classified, the information whether the 
code is optimal. In some cases we also give a reference and some further details in free text format.  

\subsection{Application programming interface}
\label{sec_api}

There is also an API available to access most data of the database. 
It is inspired by the REST (representational state transfer) style and only GET queries are supported. 
In order to access the data for the constant dimension case with parameters $q,n,d$ and $k$, you query the URL

\begin{center}
http://subspacecodes.uni-bayreuth.de/api/$q$/$n$/$d$/$k$/
\end{center}
Similarly in the mixed dimension case, the URL is
\begin{center}
http://subspacecodes.uni-bayreuth.de/api/$q$/$n$/$d$/
\end{center}
The result is a JSON file which contains a subset of the following attributes:

\begin{itemize}
\item request = contains your specified $q,n,d$ and $k$
\item \{lower,upper\}\_bound = lower or upper bound for the value $A_q(n,d;k)$
\item comments = commentaries to this entry
\item nondeduced = if the parameters are no parameters that are also viewable in the ``short'' mode, then they are trivial or computed using other parameters. nondeduced lists these other parameters.
\item \{lower,equal,upper\}\_bound\_constraints = list of tuples which contain name, parameter and value of the applied constraints
\item classified = boolean that is true if $A_q(n,d;k)$ is classified up to isomorphism
\item known codes = list of tuples of size, details, file (to enable automatic downloads) and nrisotypes (the number of isomorphism types of this entry)
\item liftedmrdsizebound = the bound for codes that contains the lifted MRD code as described in Subsection~\ref{subsec_mrd_bound}
\end{itemize}
In order to download the codes, you have to use the attribute file above and the URL
\begin{center}
http://subspacecodes.uni-bayreuth.de/codes/\emph{file}
\end{center}
We want to remark that the API (as well as the whole homepage) is still in an evolutionary phase and therefore changes are may occur.
As an example, the the interested used may open the URL
\begin{center}
\url{http://subspacecodes.uni-bayreuth.de/api/2/6/4/3/}.
\end{center}

\section{Bounds for CDCs}
\label{sec_bounds_cdcs}

For constant dimension codes much more bounds are known than for mixed dimension codes. We state lower bounds, i.e., 
constructions, in Subsection~\ref{subsec_implemented_constructions} and upper bounds in Subsection~\ref{subsec_cdc_upper}.  

\subsection{Lower bounds and constructions for CDCs}
\label{subsec_implemented_constructions}

Any subspace code is a set and hence its size is at least zero. This most trivial bound $A_q(n,d;k)\ge 0$ is \texttt{trivial\_1}. 
Lifted MRD codes, see Subsection~\ref{sec:lmrd}, are one type of building blocks of the Echelon-Ferrers construction, see 
Subsection~\ref{subsec_echelon_ferrers}. The latter is a nice interplay between the subspace distance, the rank distance and the 
Hamming distance. Another construction based on similar ideas is the so-called Coset construction, see Subsection~\ref{subsec_coset}.
The most effective general recursive construction is the so-called linkage construction and its generalization, see 
Subsection~\ref{linkage_construction}. The expurgation-augmentation method, starting from a lifted MRD code and then adding and removing 
codewords, is briefly describe in Subsection~\ref{subsubsec_expurgation}. Constant dimension codes with prescribed automorphisms are the 
topic of Subsection~\ref{subsubsec_automorphisms_lower_cdc}. Also the non-constructive lower bounds for classical codes in the 
Hamming metric can be transferred, see Subsection~\ref{subsubsec_non_constructive_lower_cdc}. Last but not least, also geometrical 
ideas can be employed in order to obtain good constructions for constant dimension codes, see Subsection~\ref{subsubsec_geometric_lower_cdc}.

\subsubsection{Lifted MRD codes} \label{sec:lmrd}

\label{subsec_lifted_MRD}
For matrices $A,B\in\mathbb{F}_q^{m\times n}$ the rank distance is defined via $d_R(A,B):=\operatorname{rk}(A-B)$. It is indeed 
a metric, as observed in \cite{gabidulin1985theory}.

\begin{Theorem}(see \cite{gabidulin1985theory})
  \label{thm_MRD_size}
  Let $m,n\ge d$ be positive integers, $q$ a prime power, and $\mathcal{C}\subseteq \mathbb{F}_q^{m\times n}$ be a rank-metric 
  code with minimum rank distance $d$. Then, $\#\mathcal{C}\le q^{\max\{n,m\}\cdot (\min\{n,m\}-d+1)}$. 
\end{Theorem}  
Codes attaining this upper bound are called maximum rank distance (MRD) codes. They exist for all (suitable) choices of parameters. If $m<d$ or $n<d$, then only $\#\mathcal{C}=1$ is possible, which may be summarized to the single upper bound 
$\#\mathcal{C}\le \left\lceil q^{\max\{n,m\}\cdot (\min\{n,m\}-d+1)}\right\rceil$. 
Using an $m\times m$ identity matrix as a prefix one obtains the so-called lifted MRD codes.

\begin{Theorem}(see \cite{silva2008rank})
  For positive integers $k,d,n$ with $k\le n$, $d\le 2\min\{k,n-k\}$, and $d\equiv 0\pmod{2}$, the size of a lifted MRD code in $\G{q}{n}{k}$ with 
  subspace distance $d$ is given by
\[
M(q,k,n,d):=q^{\max\{k,n-k\}\cdot(\min\{k,n-k\}-d/2+1)}.
\]
If $d>2\min\{k,n-k\}$, then we have $M(q,k,n,d)=1$. 
\end{Theorem}

As MRD codes can be obtained from linearized polynomials, we have the very same bound  
 implemented as \texttt{lin\_poly}:
\begin{Theorem} \label{theo:linpoly}
  (Linearized polynomials, see \cite{koetter2008coding} and Section~\ref{sec:lmrd})
\[
    A_q(n,d;k)\ge q^{(n-k)(k-d/2+1)}
\]
\end{Theorem}

\subsubsection{Echelon-Ferrers or multilevel construction}
\label{subsec_echelon_ferrers}
In \cite{etzion2009error} a generalization, the so-called multi-level construction, 
based on lifted MRD codes was presented. Let $1\le k\le n$ be integers and $b\in \mathbb{F}_2^n$ 
a binary vector of weight $k$. By $\operatorname{EF}_q(b)$ we denote the set of all $k\times n$ 
matrices over $\mathbb{F}_q$ that are in row-reduced echelon form, i.e., the Gaussian algorithm 
had been applied, and the pivot columns coincide with the positions where $b$ has a $1$-entry.

\begin{Theorem} (see \cite{etzion2009error})
  \label{thm_echelon_ferrers}
  For integers $k,n,\delta$ with $1\le k\le n$ and $1\le \delta\le \min\{k,n-k\}$, let $\mathcal{B}$ be a 
  binary constant weight code of length $n$, weight $k$, and minimum Hamming distance $2\delta$. 
  For each $b\in \mathcal{B}$ let $\mathcal{C}_b$ be a code in $\operatorname{EF}_q(b)$ 
  with minimum rank distance at least $\delta$. Then, $\cup_{b\in\mathcal{B}} \,\mathcal{C}_b$ is 
  a constant dimension code of dimension $k$ having a subspace distance of at least $2\delta$. 
\end{Theorem}  

The code $\mathcal{B}$ is also called \textit{skeleton code}. For $\mathcal{C}_b$ we have the following upper 
bound:

\begin{Theorem} (see \cite{etzion2009error})
  \label{thm_upper_bound_ef}
  Let $1\le k\le n$ be integers and $b\in \mathbb{F}_2^n$ a binary vector of weight $k$.
  Let $\mathcal{F}$ be the Ferrers diagram of $\operatorname{EF}_q(b)$ and 
  $\mathcal{C}\subseteq \operatorname{EF}_q(b)$ be a subspace code having a subspace
  distance of at least $2\delta$, then 
\[
    \#\mathcal{C} \le q^{\min\{\nu_i\,:\, 0\le i\le \delta-1\}},
\]
  where $\nu_i$ is the number of dots in $\mathcal{F}$, which are neither contained in the 
  first $i$ rows nor contained in the rightmost $\delta-1-i$ columns.  
\end{Theorem}
The authors of \cite{etzion2009error} conjecture that Theorem~\ref{thm_upper_bound_ef} is tight for 
all parameters $q$, $\mathcal{F}$, and $\delta$, which is still unrebutted. Constructions 
settling the conjecture in several cases are given in \cite{MR3480069}.

There is one rather obvious skeleton code that needs to be considered. 
Taking binary vectors with $k$ consecutive ones we are in the classical MRD case. So, taking binary vectors $v_i$, where the ones 
are located in positions $(i-1)k+1$ to $ik$ for all $1\le i\le \left\lfloor n/k\right\rfloor$, clearly gives 
a binary constant weight code of length $n$, weight $k$, and minimum Hamming distance $2k$.  

\begin{Observation} (see e.g.\ \cite{kurzspreads})
  \label{obs_multi_component} 
  For positive integers $k$, $n$ with $n>2k$ and $n\not\equiv 0\pmod{k}$, there exists a constant dimension code in $\G{q}{n}{k}$ 
  with subspace distance $2k$ having cardinality
\[
    1+\sum_{i=1}^{\left\lfloor n/k\right\rfloor-1} q^{n-ik}
    =1+q^{k+(n\,\operatorname{mod}\,k)}\cdot\frac{q^{n-k-(n\,\operatorname{mod}\,k)}-1}{q^k-1}
    =\frac{q^n-q^{k+(n\,\operatorname{mod}\,k)}+q^k-1}{q^k-1}.
\]
\end{Observation}
The observation is implemented as \texttt{multicomponent}. 
We remark that a more general construction, among similar lines and including explicit formulas for the respective cardinalities, has 
been presented in \cite{skachek2010recursive}. This lower bound for partial spreads, i.e., $d=2k$, is exactly the same as:
\begin{Theorem}
  (Partial spreads, see \cite{etzion2011error})
  If $d=2k$ then:
\[
    A_q(n,d;k)\ge \frac{q^n-q^k(q^{(n \bmod k)}-1)-1}{q^k-1}
\]
\end{Theorem}
This lower bound is implemented as \texttt{partial\_spread\_3} and equals the size of the construction of Beutelspacher, 
see \cite{beutelspacher1975partial}.

We remark that the general Echelon-Ferrers or multilevel construction contains the mentioned observation as a very easy special case. 
However, our knowledge on the size of an MRD code over $\operatorname{EF}_q(v)$ is still very limited. As mentioned, there is an 
explicit conjecture, which so far is neither proven nor disproved. Let the field size $q$, the constant dimension $k$, and the minimum 
subspace distance $d$ be fix, in order to ease the notation. By $V$ we denote the set of binary vectors of weight $k$ in $\{0,1\}^n$. 
Let $c(v)$ denote the maximum size of a known MRD code over $\operatorname{EF}_q(v)$ matching distance $d$. The optimal 
Echelon-Ferrers construction can be modeled as an ILP:
    
\begin{align*}
\max
\sum_{v\in \mathbb{F}_2^n} &c(v)\cdot x_v
\\
\operatorname{s.t.}
&x_a+x_b \le 1 &\forall a\neq b\in \mathbb{F}_2^n: d_{\operatorname{H}}(a,b)<d\\
&x_v \in \{0, 1\}                 &\forall v\in \mathbb{F}_2^n.
\end{align*}
This is implemented as \texttt{echelon\_ferrers}. However, the evaluation of this ILP is only feasible for rather moderate 
sized parameters. More sophisticated algorithmic considerations, unfortunately still unpublished, give bounds for the 
exact optimum of the Echelon-Ferrers construction, which is implemented as \texttt{ef\_computation}. A greedy-type approach 
has been considered by Alexander Shishkin, see \cite{Shishkin2014} and also \cite{shishkin2014cardinality}. It is implemented as 
\texttt{greedy\_multicomponent}. (However, we have not checked that all corresponding MRD codes for the involved Ferrers diagrams exist.) 
In \cite{gabidulin2011new,multicomponent} the authors considered block designs as skeleton codes. Further improvements on the Echelon-Ferrers 
construction can e.g.\ be found in \cite{trautmann2010new}.

By choosing some explicit skeleton code and constructing the corresponding MRD codes in $\operatorname{EF}_q(v)$, one can obtain 
explicit lower bounds:
 
\begin{Theorem}[{\cite[Example 59]{gorla2014subspace}}]
For $q > 2$: \\
$A_q(10, 6, 5) \ge q^{15} + q^{6} + 2q^{2} + q + 1$, \\
$A_q(11, 6, 5) \ge q^{18} + q^{9} + q^{6} + q^{4} + 4q^{3} + 3q^{2}$, \\
$A_q(14, 6, 4) \ge q^{20} + q^{14} + q^{10} + q^{9} + q^{8} + 2( q^{6} + q^{5} + q^{4} ) + q^{3} + q^{2}$, \\
$A_q(14, 8, 5) \ge q^{18} + q^{10} + q^{3} + 1$, and \\
$A_q(15, 10, 6) \ge q^{18} + q^{5} + 1$
\end{Theorem}
This is implemented as \texttt{Gorla\_Ravagnani\_2014}.

\medskip

The Echelon-Ferrers construction has even been fine-tuned to the so-called pending dots
\cite{etzion2013codes} implemented as \texttt{pending\_dots}, and the so-called pending 
blocks \cite{silberstein2014subspace,silberstein2013new,trautmann2013constructions} constructions. Of course, these variants have even more degrees of 
freedom, so that a general solution of the best codes within these classes of constructions is out of sight.   

Explicit series of constructions using pending dots are given by:
\begin{Theorem}[{\cite[Construction 1, see chapter IV, Theorem 16]{etzion2013codes}}]
$$
  A_q(n,2(k-1);k) \ge q^{2(n-k)} + A_q(n-k,2(k-2),k-1)
$$ 
if $q^2+q+1 \ge s$ with $s=n-4$ if n is odd and $s=n-3$ else
\end{Theorem}
This is implemented as \texttt{construction\_1}.

\begin{Theorem}[{\cite[Construction 2, see chapter IV, Theorem 17]{etzion2013codes}}]
$$
  A_q(n,4;3) \ge q^{2(n-3)} + \sum_{i=1}^{\alpha} q^{2(n-3-(q^2+q+2)i)}
$$ 
if $q^2+q+1 < s$ with $s=n-4$ if n is odd and $s=n-3$ else and $\alpha = \left\lfloor \frac{n-3}{q^2+q+2} \right\rfloor$
\end{Theorem}
This is implemented as \texttt{construction\_2}.

\medskip

Explicit series of constructions using pending blocks are given by:
\begin{Theorem}[{\cite[Construction A, see chapter III, Theorem 19, Corollary 20]{silberstein2014subspace}}]
Let $n\geq \frac{k^2+3k-2}{2}$ and $q^2+q+1\geq \ell$, where $\ell= n-\frac{k^2+k-6}{2}$ for odd $n-\frac{k^2+k-6}{2}$ (or $\ell= n-\frac{k^2+k-4}{2}$ for even $n-\frac{k^2+k-6}{2}$). Then
$A_q(n, 2k-2; k) \ge q^{2(n-k)}+\sum_{j=3}^{k-1} q^{2(n-\sum_{i=j}^k i)}+\gauss{n-\frac{k^2+k-6}{2}}{2}{q}$.
\end{Theorem}
This is implemented as \texttt{construction\_ST\_A\_1}.

\begin{Theorem}[{\cite[Construction B, see chapter IV, theorem 26, Corollary 27]{silberstein2014subspace}}]
Let $n\geq 2k+2$. Then
$A_q(n,4;k) \ge \sum_{i=1}^{\lfloor\frac{n-2}{k}\rfloor -1}\left( q^{(k-1)(n-ik)}+ \frac{(q^{2(k-2)}-1)(q^{2(n-ik-1)}-1)}{(q^4-1)^2}q^{(k-3)(n-ik-2)+4}\right)$.
\end{Theorem}
This is implemented as \texttt{construction\_ST\_B}.

\begin{Theorem}
Let $q \ge 2$ be a prime power and $2 \le d/2 \le k \le v-k$ integers.
If additionally $d \le k+1$, then
\[A_q(n,d;k)\ge q^{(v-k)(k-d/2+1)} \frac{q^{(d/2)^2(M+1)}-1}{q^{(d/2)^2}-1}q^{-(d/2)^2M}\]
with $M=\lceil 2(v-k)/d \rceil$.
\end{Theorem}

This is implemented as \texttt{two\_pivot\_block\_construction}.

\subsubsection{Coset construction}
\label{subsec_coset}
The so-called Coset construction, see \cite{heinlein2015coset}, grounds, similar as the Echelon-Ferrers construction, on the interplay between  
the subspace distance, the rank distance and the Hamming distance. Another way to look at it is that it generalizes the 
construction from \cite[Theorem~18]{etzion2013codes}, which yields $A_2(8,4;4)\ge 4797$. 
Implemented as \texttt{construction\_3}, we have for general prime powers $q$: 
\begin{Theorem}[{\cite[Construction 3, see chapter V, Theorem 18, Remark 6]{etzion2013codes}}]
$A_q(8,4;4) \ge q^{12}+\gauss{4}{2}{q}(q^2+1)q^2+1$.
\end{Theorem}

Even more than the Echelon-Ferrers construction, it 
is a rather general approach that restricts the general optimization problem of determining the subspace codes with the maximum cardinality 
to the best combination of some structured building blocks. Here the building blocks are even more sophisticated than the 
MRD codes over $\operatorname{EF}_q(v)$, so that in general only lower and upper bounds for their sizes are known. Two explicit 
parameterized constructions are given by;

\begin{Theorem}[{\cite[Section V-A, Theorem 11]{heinlein2015coset}}]
For all $q$, we have $A_q(8,4;4)\ge q^{12}+\gauss{4}{2}{q}(q^2+1)q^2+1$.

For each $k\ge 4$ and arbitrary $q$ we have $A_q(3k-3,2k-2;k)\ge q^{4k-6}+\frac{q^{2k-3}-q}{q^{k-2}-1}-q+1$.
\end{Theorem}

\begin{Theorem}
We have $A_2(18,6;9) \ge 9241456945250010249$.
\end{Theorem}

These theorems are implemented as \texttt{coset\_construction}.

\begin{Theorem}[{\cite[Theorem~9]{heinlein2015coset}}]
If $\gauss{\mathbb{F}_q^{n_i}}{k_i}{}$ admit parallelisms, i.e., a partition into spreads, for $i=1,2$ then $A_q(n_1+n_2,4;k_1+k_2) \ge s_1 \cdot s_2 \cdot \min\{p_1,p_2\} \cdot m$, where $s_i=\frac{q^{n_i}-1}{q^{k_i}-1}$ is the size of a spread and $p_i=\frac{\gauss{n_i}{k_i}{q}}{s_i}$ is the size of a parallelism in $\gauss{\mathbb{F}_q^{n_i}}{k_i}{}$ for $i=1,2$, and $m=\lceil q^{\max\{k_1,n_2-k_2\}(\min\{k_1,n_2-k_2\}-1)} \rceil$ is the size of an MRD code with shape $k_1 \times (n_2-k_2)$ and rank distance $2$ over $\mathbb{F}_q$.
\end{Theorem}

Unfortunately, the existence question for parallelisms in $\gauss{\mathbb{F}_q^n}{k}{}$ is still open in general. They are known to exist for:
\begin{enumerate}
\item $q=2, k=2$ and $n$ even;
\item $k=2$, all $q$ and $n=2^m$ for $m\ge 2$;
\item $n=4$, $k=2$, and $q\equiv 2\pmod 3$;
\item $q=2, k=3, n=6$,
\end{enumerate}
see e.g.\ \cite{etzion2015galois}. All applicable parameter combinations for $(n_1,k_1)$ and $(n_2,k_2)$ are implemented as 
\texttt{coset\_construction\_parallelism\_part}.

Based on a packing of a $(6,77,4;3)$ code into several subcodes with minimum subspace distance $6$, a construction 
for $A_2(10,6;4)\ge 4173$ was obtained in \cite[Theorem 13]{heinlein2015coset}. This is still the best known code for 
these parameters and can be downloaded as a file. For $q\ge 3$ it remains unknown whether a similar construction can improve 
upon the best known construction obtained from the Echelon-Ferrers construction. 

\subsubsection{Linkage constructions}
\label{linkage_construction}

A powerful construction to obtain large codes from a given code $\mathcal{C}$ is to append all possible choices of an MRD code to 
the matrices in row-echelon form of the codewords of $\mathcal{C}$. This resulting size of the constructed code is the size of $\mathcal{C}$ 
times the size of the MRD code. This approach is called Construction~D in \cite{silberstein2014subspace}, 
see \cite[Theorem 37]{silberstein2014subspace} and also \cite[Theorem 5.1]{gluesing2015cyclic} and is implemented as \texttt{construction\_D}. 

Performing a tighter analysis of the occurring subspace distances one notices that one can add further codewords from a code 
in a smaller ambient space to Construction D. This gives:

\begin{Theorem}
(linkage by Silberstein and (Horlemann-)Trautmann, see \cite[Corollary~39]{silberstein2014subspace})
For $3k \le n$ and $k \le \Delta \le n$ we have:
\[A_q(n,d;k) \ge q^{\Delta(k-d/2+1)}A_q(n-\Delta,d;k)+A_q(\Delta,d;k)\]
\end{Theorem}
This bound is implemented as \texttt{linkage\_ST}. 
Without the assumption $3k\le n$, the same bound is independently obtained in:
\begin{Theorem}
(linkage by Gluesing-Luerssen and Troha \cite[Theorem~2.3]{MR3543532})
For $k \le m \le n-k$ we have:
\[A_q(n,d;k) \ge A_q(m,d;k) \cdot \left\lceil q^{ (n-m)(k-d/2+1) } \right\rceil + A_q(n-m,d;k)\]
\end{Theorem}
This bound is implemented as \texttt{linkage\_GLT}. We remark that for $n<3k$ better constructions are known, see e.g.\ \cite[Footnote 2]{heinlein2017asymptotic}. 
An improved analysis of the involved distances yields: 
\begin{Theorem} (improved linkage)\label{theo:improved_linkage}
For $k \le m \le n-d/2$ we have:
\[A_q(n,d;k) \ge A_q(m,d;k) \cdot \left\lceil q^{ \max\{n-m,k\}(\min\{n-m,k\}-d/2+1) } \right\rceil + A_q(n-m+k-d/2,d;k)\]
\end{Theorem}
This bound is implemented as \texttt{improved\_linkage}.
The description of the application of all three constraints contains $\Delta$ respective $m$ in brackets.

Further refinements of the linkage construction are described in \cite{kurz2019note}, se also \cite{heinlein2019generalized}. 
\begin{Theorem}(\cite[Proposition 4.6]{kurz2019note})
  \begin{eqnarray*}
    A_q(12,4;4)&\ge & q^{24}+q^{20}+q^{19}+3q^{18}+2q^{17}+3q^{16}+q^{15}+q^{14}+q^{12}\\&&+q^{10}+2q^8+2q^6+2q^4+q^2\\ 
    A_q(13,4;4) &\ge& q^{27}+q^{23}+q^{22}+3q^{21}+2q^{20}+3q^{19}+q^{18}+q^{17}+q^{15}\\&&+q^{12}+q^{10}+q^9+q^8+q^7+q^6+q^5+q^3
  \end{eqnarray*}
\end{Theorem}
This bound is implemented as \texttt{noted\_linkage\_special}.

\begin{Theorem}(\cite[Theorem 4.2]{kurz2019note})
  We have
  $$
    A_q(n,d;k)\ge A_q(m,d;k)\cdot q^{2(n-m)}+A_q(n-m,d-2;k-1)
  $$
  for $d=2k-2$ and $3\le k\le m\le n-k$. 
\end{Theorem}
This bound is implemented as \texttt{noted\_linkage}. Note that the description contains the parameter $m$ in brackets.

\cite[Theorem 21]{heinlein2019generalized} is implemented as \texttt{generalized\_linkage}.

\subsubsection{The expurgation-augmentation method}
\label{subsubsec_expurgation}
The success of the Echelon--Ferrers and the coset construction is mainly given by the fact that lifted MRD codes 
have a quite large cardinality, which is asymptotically optimal in a certain sense. While both two methods try to append 
some additional subcodes, the linkage constructions employ the MRD codes in a product type construction. Another approach 
is to start from a lifted MRD code, remove some codewords in order to add more codewords again. This approach is coined 
\textit{expurgation-augmentation} and invented by Thomas Honold.

The starting point is possible a computer--free construction for the lower bound $A_2(7,4;3)\ge 329$, see \cite{haiteng2014poster}, 
which was previously obtained by a computer search using prescribed automorphisms, see \cite{braun2014q}. Successors are:

\begin{Theorem}[{\cite[Theorem 2]{hkk77}}] \label{theo:HonoldKiermaierKurz_n6_d4_k3}
$A_q(6,4;3) \ge q^6+2q^2+2q+1 $ for $3 \le q$.
\end{Theorem}
This is implemented as \texttt{HonoldKiermaierKurz\_n6\_d4\_k3}.
Note that the right hand side in Theorem~\ref{theo:HonoldKiermaierKurz_n6_d4_k3} is larger for all $q \ge 3$ than the right hand 
side in Theorem~\ref{theo:CossidentePavese_n6_d4_k3}.

\begin{Theorem}[\cite{honold15talkalcoma}]
$A_q(7,4;3) \ge q^8 + q^5 + q^4 - q - 1$
\end{Theorem}
This is implemented as \texttt{construction\_honold} and superseded by \texttt{construction\_HK15}.

\begin{Theorem}[{\cite[Theorem 4]{MR3444245}}]
$A_2(7,4;3) \ge 329$, $A_3(7,4;3) \ge 6977$, $A_q(7,4;3) \ge q^{8} + q^{5} + q^{4} + q^{2} - q$
\end{Theorem}
This is implemented as \texttt{construction\_HK15}.

\medskip

While the sketched idea of the expurgation-augmentation method is rather general, several theoretical 
insights are possible. Prescribing automorphisms in the constructions also helps to obtain 
optimization problems that are more structured and computationally feasible. A whole theoretical 
framework is introduced in \cite{ai2016expurgation}. As a purely analytical result we have: 

\begin{Theorem}[{\cite[Main Theorem]{ai2016expurgation}}]
$A_2(n,4;3) \geq 2^{2(n-3)}+\frac{9}{8}\gauss{n-3}{2}{2}$ for $n \equiv 7\pmod{8}$ \\
$A_2(n,4;3)\geq 2^{2(n-3)}+\frac{81}{64}\gauss{n-3}{2}{2}$ for $n\equiv 3\pmod{8}$ and $n\geq 11$.
\end{Theorem}
This is implemented as \texttt{expurgation\_augmentation\_general}.

Explicit computer calculations allow further improvements:
\begin{Theorem}[{\cite[Table 1]{ai2016expurgation}}]
$A_2(7,4;3) \ge 2^{8} + 45$, \\
$A_2(8,4;3) \ge 2^{10} + 93$, \\
$A_2(9,4;3) \ge 2^{12} + 756$, \\
$A_2(10,4;3) \ge 2^{14} + 2540$, \\
$A_2(11,4;3) \ge 2^{16} + 13770$, \\
$A_2(12,4;3) \ge 2^{18} + 47523$, \\
$A_2(13,4;3) \ge 2^{20} + 239382$, \\
$A_2(14,4;3) \ge 2^{22} + 775813$, \\
$A_2(15,4;3) \ge 2^{24} + 3783708$, and \\
$A_2(16,4;3) \ge 2^{26} + 12499466$
\end{Theorem}
This is implemented as \texttt{expurgation\_augmentation\_special\_cases}.

\subsubsection{Codes with prescribed automorphisms}
\label{subsubsec_automorphisms_lower_cdc}

The computational complexity of the general optimization problem for $A_q(n,d;k)$ can be reduced if one assumes that 
the desired constant dimension code $\mathcal{C}$ admits some automorphisms, see \cite{kohnert2008construction}.
So, the idea is to prescribe some subgroup $G$ of the automorphism group. If $G$ is cyclic, then some authors 
speak of cyclic orbit codes, see e.g.\ \cite{climent2017construction,gluesingcardinality,gluesing2015cyclic,horlemann2017correction,
trautmann2010orbit,trautmann2013cyclic,rosenthal2013complete,manganiello2011conjugacy,trautmann2011complete,gluesing2019distance}. 
For these objects one can utilize the theory of subspace polynomials, see \cite{ben2016subspace,otal2017cyclic}, and Sidon 
spaces, see \cite{roth2017construction}. The Singer cycle is one prominent 
example since it acts transitively on the one-dimensional subspaces of $\mathbb{F}_q^n$. We restate the 
computational results from~\cite{kohnert2008construction} for $A_2(n,4;3)$: 

\begin{center}\begin{tabular}{|c|c|c|c|c|c|c|}
\hline 
$n$&
$k$&
$l$&
$\#$ orbits&
$\#$ codewords&
$d$\tabularnewline
\hline
\hline 
$6$&
$3$&
$1$&
$19$&
$1\cdot63=63$&
$4$\tabularnewline
\hline 
$7$&
$3$&
$2$&
$93$&
$2\cdot127=254$&
$4$\tabularnewline
\hline 
$8$&
$3$&
$5$&
$381$&
$5\cdot255=1275^{*}$&
$4$\tabularnewline
\hline 
$9$&
$3$&
$11$&
$1542$&
$11\cdot511=5621^{*}$&
$4$\tabularnewline
\hline 
$10$&
$3$&
$21$&
$6205$&
$21\cdot1023=21483^{*}$&
$4$\tabularnewline
\hline
$11$&
$3$&
$39$&
$24893$&
$39\cdot2047=79833^{*}$&
$4$\tabularnewline
\hline
$12$&
$3$&
$77$&
$99718$&
$77\cdot4095=315315^{*}$&
$4$\tabularnewline
\hline
$13$&
$3$&
$141$&
$399165$&
$141\cdot8191=1154931$&
$4$\tabularnewline
\hline
$14$&
$3$&
$255$&
$1597245$&
$255\cdot16383=4177665$&
$4$\tabularnewline
\hline
\end{tabular}\par\end{center}
 
Here $l$ denotes the number of chosen orbits from the total number of orbits. Those code sizes that were the best known 
lower bound at that time are marked with an asterisk. We remark that the stated values correspond to the optimal solutions 
of the corresponding ILP for $6\le n\le 8$. For $n=9$ it was reported that $l=12$ might be possible, which would be larger 
than the best known code for $A_2(9,4;3)\ge 5986$ found in \cite{new_lower_bounds_cdc}. Later on the following, partially weaker,  
results have been obtained using the normalizer of the Singer cycle:
\begin{Theorem}[{\cite[Example 2.7 and 2.8]{MR3431963}}]
$A_2(n,4;3) \ge n \cdot (2^n-1)$ for $12 \le n \le 20$, \\
$A_2(8,4;3) \ge 2 \cdot (2^8-1)$, \\
$A_2(9,4;3) \ge 9 \cdot (2^9-1)$, \\
$A_2(13,6;4) \ge 13 \cdot (2^{13}-1)$, and \\
$A_2(17,6;4) \ge 17 \cdot (2^{17}-1)$
\end{Theorem}
This is implemented as \texttt{Bardestani\_Iranmanesh}.

\medskip

A slight variation of cyclic subspace codes was considered in \cite{garcia2016quasi,gutierrez2015some}.

\subsubsection{Transferred \emph{classical} non-constructive lower bounds}
\label{subsubsec_non_constructive_lower_cdc}

The classical Gilbert-Varshamov lower bound, based on sphere coverings, has been transferred to constant dimension codes:
\begin{Theorem}
  \label{thm_cdc_sphere_covering}
  (Sphere covering, see \cite{koetter2008coding})
\[
    A_q(n,d;k)\ge \gauss{n}{k}{q}/\left(\sum_{i=0}^{(d/2-1)+1} \gauss{k}{i}{q} \cdot \gauss{n-k}{i}{q} \cdot q^{i^2}\right)
\]
\end{Theorem}
This lower bound is implemented as \texttt{sphere\_covering}.

\medskip

A Graham-Sloane type bound was obtained in \cite{xia2008graham}:
\begin{Theorem}
  (Graham, Sloane, see \cite{xia2008graham})
\[
    A_q(n,d;k)\ge \frac{(q-1)\gauss{n}{k}{q}}{(q^n-1)q^{n(d/2-2)}}
\]
\end{Theorem}
This lower bound is implemented as \texttt{graham\_sloane}. For minimum subspace distance $d=4$ is yields a 
strictly larger lower bound than Theorem~\ref{thm_cdc_sphere_covering}.  

\subsubsection{Geometric constructions}
\label{subsubsec_geometric_lower_cdc}

Geometric concepts like the Segre variety and the Veronese variety where also used to obtain constructions for 
constant dimension codes:

\begin{Theorem}[{\cite[Theorem 3.11]{cossidente2016subspace}}]
If $n \ge 5$ is odd, then

$ A_q(2n,4;n) \ge q^{n^2-n} + \sum_{r=2}^{n-2} \gauss{n}{r}{q} \sum_{j=2}^{r} (-1)^{(r-j)} \gauss{r}{j}{q} q^{\binom{r-j}{2}}(q^{n(j-1)}-1) + \prod_{i=1}^{n-1} (q^i+1) - q^{\frac{n(n-1)}{2}} - \gauss{n}{1}{q} \left( q^{\frac{(n-1)(n-2)}{2}} - q^{\frac{(n-1)(n-3)}{4}} \prod_{i=1}^{\frac{n-1}{2}} (q^{2i-1}-1) \right) +y(y-1) + 1$, using $y:=q^{n-2}+q^{n-4}+\dots+q^3+1$.
\end{Theorem}
This is implemented as \texttt{CossidentePavese14\_theorem311}.

\begin{Theorem}[{\cite[Theorem 3.8]{cossidente2016subspace}}]
If $n \ge 4$ is even, then
\begin{align*}
A_q(2n,4;n) \ge q^{n^2-n} + \sum_{r=2}^{n-2} \gauss{n}{r}{q} \sum_{j=2}^{r} (-1)^{(r-j)} \gauss{r}{j}{q} q^{\binom{r-j}{2}}(q^{n(j-1)}-1)
\\+ (q+1) \left( \prod_{i=1}^{n-1} (q^i+1) - 2q^{\frac{n(n-1)}{2}} + q^{\frac{n(n-2)}{4}} \prod_{i=1}^{\frac{n}{2}} (q^{2i-1}-1) \right) - q \cdot |G| + \gauss{\frac{n}{2}}{1}{q^2} \left( \gauss{\frac{n}{2}}{1}{q^2} - 1 \right) + 1
\end{align*}
using $|G| = 2 \prod_{i=1}^{n/2-1}(q^{2i}+1) - 2q^{(n(n-2)/4)} $ if $n/2$ is odd and $|G| = 2 \prod_{i=1}^{n/2-1}(q^{2i}+1) - 2q^{(n(n-2)/4)} + q^{n(n-4)/8}\prod_{i=1}^{n/4}(q^{4i-2}-1) $ if $n/2$ is even.
\end{Theorem}
This is implemented as \texttt{CossidentePavese14\_theorem38}.

\begin{Theorem}[{\cite[Theorem 4.3]{cossidente2016subspace}}]
$A_q(8,4;4) \ge q^{12}+q^2(q^2+1)^2(q^2+q+1)+1$
\end{Theorem}
This is implemented as \texttt{CossidentePavese14\_theorem43}.

\begin{Theorem}[{\cite[Corollary 7.4]{MR3544049}}] \label{theo:CossidentePavese_n6_d4_k3}
$A_q(6,4;3) \ge q^3(q^2-1)(q-1)/3 + (q^2+1)(q^2+q+1)$
\end{Theorem}
This is implemented as \texttt{CossidentePavese\_n6\_d4\_k3}.

A survey on geometric methods for constant dimension codes can be found in \cite{cossidente2018geometrical}.

\subsubsection{Further lower bounds}

\begin{Theorem}[{\cite[Lower bound of Lemma 2.4]{MR3543542}}] 
\(1\le d/2\leq k\le\lfloor n/2\rfloor, d\equiv 0\bmod 2 \Rightarrow A_q(n,d;k)>q\cdot A_q(n,d;k-1)\).
\end{Theorem}

This is implemented as \texttt{HKK\_theorem\_3\_3\_i\_lower\_bound\_A\_q(6,4)}.

\begin{Theorem}(\cite[Theorem 3.13]{cossidente2019subspace})
  $$
    A_q(9,4;3)\ge q^{12}+2q^8+2q^7+q^6+q^5+q^4+1
  $$
\end{Theorem}
This bound is implemented as \texttt{CossidenteMarinoPavese2019\_T313}.

\begin{Theorem}(\cite[Theorem 4.2 and 4.7]{cossidente2019subspace})
  $$
    A_q(6,4;3)\ge \left(q^3-1\right)\left(q^2+q+1\right)
  $$
\end{Theorem}
This bound is implemented as \texttt{CossidenteMarinoPavese2019\_T42\_T47}.

\begin{Theorem}(\cite[Corollary 4]{kurz2019subspaces})
  $$
    A_q(9,4;3)\ge q^{12}+2q^8+2q^7+q^6+2q^5+2q^4-2q^2-2q+1
  $$
\end{Theorem}
This bound is implemented as \texttt{Kurz20192\_C4}.

Further lower bounds, that are not implemented yet, can be found in \cite{kurz2019subspaces,cossidente2019combining,
chen2019new,cossidente2016subspaceDCC,he2019construction,he2019hierarchical,liu2019parallel,antrobus2019maximal,
liu2019constructions,liu2019several,zhang2019constructions}.

\subsection{Upper bounds for CDCs}
\label{subsec_cdc_upper}

\noindent
Surveys and partial comparisons of upper bounds for constant dimension codes can e.g.\ be found in  
\cite{MR3063504,heinlein2017asymptotic,khaleghi2009subspace}.

Assuming $0\le k\le n$ we always have $A_q(n,d;k)\ge 1$.
Since we can take no more than all subspaces of a given dimension, 
we obtain the trivial upper bound $A_q(n,d;k)\le\gauss{n}{k}{q}$ which is implemented as \texttt{all\_subs}. 
Transferred bounds from classical coding theory are stated in Subsection~\ref{subsec_classical_upper_bounds}. Of 
special importance is the Johnson bound, so that implications are treated in Subsection~\ref{subsubsec_Johnson}. 
In our description of known constructions we have seen that the lifted MRD codes play a major role in many constructions.  
For those codes tighter upper bounds are known, see Subsection~\ref{subsec_mrd_bound}. As the Johnson bound recurs back 
to bounds for partial spreads we state the corresponding bounds in Subsection~\ref{subsec_upper_bounds_spreads}. Everything 
else that does not fit into the previous categories is collected in Subsection~\ref{subsec_further_upper}. 
 
Nevertheless there is a large variety of upper bounds for constant dimension codes, the picture for the currently tightest known 
bounds is pretty clear. Besides the exact values $A_2(6,4;3)=77$ \cite{hkk77} and $A_2(8,6;4)=257$ \cite{heinlein2017classifying}, obtained with integer linear programming 
techniques, see Subsection~\ref{subsec_further_upper}, all upper bounds are given by formula (\ref{ie_best_upper_bound}), which refers 
to partial spreads. For partial spreads, Theorem~\ref{thm_spread} (the construction for spreads), Theorem~\ref{ps_two_mac_williams}, 
Theorem~\ref{ps_three_mac_williams}, and Theorem~\ref{ps_four_mac_williams} are sufficient. The latter three results are implications 
of the Delsarte linear programming method for projective linear divisible codes with respect to the Hamming metric. See 
\cite{honold2016partial} for further details.  
 
\subsubsection{\textit{Classical} coding theory bounds}
\label{subsec_classical_upper_bounds}

\begin{Theorem}(Singleton bound, see \cite{koetter2008coding})
\[
    A_q(n,d;k)\le \gauss{n-d/2+1}{k-d/2+1}{q}
\]
\end{Theorem}
This upper bound is implemented as \texttt{singleton}.

\begin{Theorem}(Sphere packing bound, see \cite{koetter2008coding})
\label{thm:sphere_packing}
\[
    A_q(n,d;k)\le \left\lfloor \gauss{n}{k}{q}/\left(\sum_{i=0}^{\left\lfloor(d/2-1)/2\right\rfloor} \gauss{k}{i}{q} \cdot \gauss{n-k}{i}{q} \cdot q^{i^2}\right) \right\rfloor
\]
\end{Theorem}
This upper bound is implemented as \texttt{sphere\_packing}.

\begin{Theorem}(Anticode bound, see \cite{etzion2011error})\label{theo:anticode}
\[
    A_q(n,d;k)\le \left\lfloor \gauss{n}{k}{q}/\gauss{n-k+d/2-1}{d/2-1}{q} = \frac{\gauss{n}{k-d/2+1}{q}}{\gauss{k}{k-d/2+1}{q}}\right\rfloor
\]
\end{Theorem}
This upper bound is implemented as \texttt{anticode}.

\begin{Theorem}[{\cite[Proposition 3]{zhang2011linear}}]
\label{lin_prog_bound_cdc}
For integers \(0 \le k \le n\) and \(2 \le d \le \min\{k,n-k\}\) such that \(d\) is even, we have  \begin{equation}
    A_q(n,d;k)\le \max \left.\left\{ 1+\sum_{i=d/2}^k x_i \,\right|\, \sum_{i=d/2}^k -Q_j(i) x_i \le u_j \,\forall j=1, 2, \ldots, k \text{ and } x_i \ge 0 \,\forall i=d / 2, d/2+1, \ldots, k \right\}
\end{equation}
 with
\begin{equation}
 u_j=\gauss{n}{j}{q}-\gauss{n}{j-1}{q},
\end{equation}
\begin{equation}
 v_i=q^{i^2}\gauss{l}{i}{q}-\gauss{n-1}{i}{q},
\end{equation}
\begin{equation}
 E_i(j)=\sum_{m=0}^i (-1)^{i-m} q^{\binom{i-m}{2}+jm}\gauss{k-m}{k-1}{q}\gauss{k-j}{m}{q}\gauss{n-k-j+m}{m}{q}\text{ and}
\end{equation}
\begin{equation}
 Q_j(i)=\frac{u_j}{v_i}E_i(j)
\end{equation}
\end{Theorem} 
 This is implemented as \texttt{linear\_programming\_bound}.

\medskip

In 1962 Johnson obtained several bounds for constant weight codes, see~\cite{johnson1962new}. All of them could be 
transferred to constant dimension codes:

\begin{Theorem}[{\cite[Theorem 2]{xia2009johnson}}]
\label{thm_johnson_I}
$A_q(n,d;k) \le \left\lfloor \frac{(q^k-q^{k-d/2})(q^n-1)}{(q^k-1)^2-(q^n-1)(q^{k-d/2}-1)} \right\rfloor$ if $(q^k-1)^2-(q^n-1)(q^{k-d/2}-1) > 0$
\end{Theorem}
This is implemented as \texttt{XiaFuJohnson1}. We remark that requested condition can be simplified considerably, see \cite{heinlein2017asymptotic}  
and also \cite{multicomponentgabidulin2016}.
\begin{Proposition}
  \label{prop_johnson_I}
  For $0\le k<n$, the bound in Theorem~\ref{thm_johnson_I} is applicable iff $d=2\min\{k,n-k\}$ and $k\ge 1$. Then, it is equivalent to
  \[
    A_q(n,d;k)\le \frac{q^v-1}{q^{\min\{k,n-k\}}-1}.
  \]
\end{Proposition}
In other words, Theorem~\ref{thm_johnson_I} is a, rather weak, bound for partial spreads obtained by dividing the number of points of the 
ambient space by the number of points of the codewords.

\begin{Theorem}(Johnson bounds, see \cite{etzion2011error})
\label{thm_johnson_cdc}
\begin{equation}
  \label{ie_j_2}
  A_q(n,d;k)\le \left\lfloor \frac{(q^n-1)\cdot A_q(n-1,d;k-1)}{q^k-1}\right\rfloor
\end{equation}
 and
\begin{equation}
  \label{ie_j_o}
  A_q(n,d;k)\le \left\lfloor \frac{(q^n-1)\cdot A_q(n-1,d;k)}{q^{n-k}-1}\right\rfloor
\end{equation}
\end{Theorem}
These upper bounds are implemented as \texttt{johnson\_1} and \texttt{johnson\_2}, respectively. Note that for $d=2k$ 
Inequality~(\ref{ie_j_2}) gives $A_q(n,2k;k)\le \left\lfloor\frac{q^n-1}{q^k-1}\right\rfloor$ since we have $A_q(n-1,2k;k-1)=1$ 
by definition. Similarly, for $d=2(n-k)$, Inequality~(\ref{ie_j_o}) gives $A_q(n,2n-2k;k)\le \left\lfloor\frac{q^n-1}{q^{n-k}-1}\right\rfloor$. 
Some sources like \cite[Theorem~3]{xia2009johnson} list just Inequality~\ref{ie_j_2} and omit 
Inequality~\ref{ie_j_o}. This goes in line with the treatment of the classical Johnson type bound II for binary error-correcting codes, 
see e.g.\ \cite[Theorem~4 on page 527]{MR0465510}, where the other bound is formulated as Problem~(2) on page 528 with the hint that ones should 
be replaced by zeros. Analogously, we can consider orthogonal codes:
\begin{Proposition}[{\cite[Proposition 2]{heinlein2017asymptotic}, cf.~\cite[Section III, esp.~Lemma~13]{etzion2011error}}]
  Inequality~(\ref{ie_j_2}) and Inequality~(\ref{ie_j_o}) are equivalent using orthogonality.
\end{Proposition}   

\subsubsection{Implications and generalizations of the Johnson bounds}
\label{subsubsec_Johnson}

The constraints of the binary linear program

\begin{align*}
\max
\sum_{U \in \gauss{\mathbb{F}_q^n}{k}{\phantom{q}}} &x_U
\\
\operatorname{s.t.}
\sum_{U \ge W} &x_U \le A_q(n-w,d;k-w) &\forall W \in \gauss{\mathbb{F}_q^n}{w}{\phantom{q}} \forall w \in \{1, \ldots, k-1\}
\\
\sum_{U \le A} &x_U \le A_q(a,d;k)   &\forall A \in \gauss{\mathbb{F}_q^n}{a}{\phantom{q}} \forall a \in \{k+1, \ldots, n-1\}
\\
&x_U \in \{0, 1\}                 &\forall U \in \gauss{\mathbb{F}_q^n}{k}{\phantom{q}}
\end{align*}

can be combined to get:

\begin{tabular}{ll}
$A_q(n,d;k) \le \frac{\gauss{n}{w}{q}}{\gauss{k}{w}{q}}A_q(n-w,d;k-w) \;\forall w \in \{1,\ldots,k-d/2\}$     & \texttt{ilp\_1}    \\
$A_q(n,d;k) \le \frac{\gauss{n}{w}{q}}{\gauss{k}{w}{q}} \;\forall w \in \{k-d/2+1, \ldots, k-1\}$     & \texttt{ilp\_2}    \\
$A_q(n,d;k) \le \frac{\gauss{n}{a}{q}}{\gauss{n-k}{a-k}{q}} \;\forall a \in \{k+1, \ldots, k+d/2-1\}$     & \texttt{ilp\_3}    \\
$A_q(n,d;k) \le \frac{\gauss{n}{a}{q}}{\gauss{n-k}{a-k}{q}}A_q(a,d;k) \;\forall a \in \{k+d/2,\ldots,n-1\}$     & \texttt{ilp\_4}    \\
\end{tabular}

Note that \texttt{ilp\_2} is \texttt{ilp\_1} using $A_q(n-w,d;k-w)=1$ for $w \in \{k-d/2+1, \ldots, k-1\}$. The same is true for \texttt{ilp\_3} 
and \texttt{ilp\_4} using $A_q(a,d;k)=1$ for $a \in \{k+1, \ldots, k+d/2-1\}$. 
Also, note that \texttt{ilp\_1} is for $w=1$ \texttt{johnson\_1}, \texttt{ilp\_2} is for $w=k-d/2+1$ \texttt{anticode}, and \texttt{ilp\_4} is 
for $a=n-1$ \texttt{johnson\_2}. In general, all these upper bounds are obtained from iterative applications of the Johnson bound from 
Theorem~\ref{thm_johnson_cdc}. As it turns out that this bound is one of the tightest known bounds, we look at it in more detail.    
In the classical Johnson space the optimal combination of the corresponding two inequalities in a recursive application is unclear,  
see e.g.\ \cite[Research Problem 17.1]{MR0465510}. For constant dimension codes there is an easy criterion for the optimal choice: 
\begin{Proposition}[{\cite[Proposition 3]{heinlein2017asymptotic}}]
  \label{prop_optimal_johnson}
  For $k \le n/2$ we have
  \[
    \left\lfloor \frac{q^n-1}{q^k-1} A_q(n-1,d;k-1) \right\rfloor
    \le
    \left\lfloor \frac{q^n-1}{q^{n-k}-1} A_q(n-1,d;k) \right\rfloor,
  \]
  where equality holds iff $n=2k$.
\end{Proposition}

With this the following non-recursive upper bound can be obtained:
\begin{equation}
  \label{ie_johnson_to_partial_spread}
  A_q(n,d;k)\le \left\lfloor \frac{q^n-1}{q^k-1}\cdot\left\lfloor \frac{q^{n-1}-1}{q^{k-1}-1}\cdot
  \left\lfloor \dots \cdot\left\lfloor\frac{q^{n'+1}-1}{q^{d'+1}-1}\cdot A_q(n',d;d') 
  \right\rfloor \dots \right\rfloor \right\rfloor\right\rfloor,
\end{equation}
where $d'=d/2$ and $n'=n-k+d'$, i.e., the unknown value on the right hand side corresponds to the case of a partial spread. Some authors 
plug in Theorem~\ref{spread_bound} in order to obtain an explicit upper bound. However, for partial spreads tighter bounds are 
available, see Subsection~\ref{subsec_upper_bounds_spreads}.

There is a recent improvement of Inequality~(\ref{ie_johnson_to_partial_spread}) or Theorem~\ref{thm_johnson_cdc}. The idea behind 
Inequality~(\ref{ie_j_2}) is that we can recursively determine an upper bound $\Lambda$ on the number of codewords that can be incident with a  
given point. With this the number of codewords is at most $\Lambda\gauss{n}{1}{q}/\gauss{k}{1}{q}$. If this number is not an integer 
it can be rounded down. In that case it means that some points are not incident to $\Lambda$ codewords. Consider the multiset of points 
with multiplicity $\Lambda$ minus the number of incidences with codewords. This multiset is equivalent to a linear 
code over $\mathbb{F}_q$, whose Hamming weights are divisible by $q^{k-1}$, see \cite{kiermaier2017improvement}. Actually, this 
is a generalization of the concept of holes and linear projective divisible codes, see Subsection~\ref{subsec_upper_bounds_spreads}.    
\begin{Theorem}[{\cite[Theorem~3 and Theorem~4]{kiermaier2017improvement}}]
Let 
\[
m= \gauss{n}{1}{q}\cdot A_q(n-1,d;k-1)-\gauss{k}{1}{q}\cdot \left\lfloor \frac{\gauss{n}{1}{q}\cdot A_q(n-1,d;k-1)}{\gauss{k}{1}{q}}\right\rfloor + 
\gauss{k}{1}{q} \cdot \delta
\]
for some $\delta\in\mathbb{N}_0$. If no $q^{k-1}$-divisible multiset of points in $\mathbb{F}_q^v$ of cardinality $m$ exists, then 
\[A_q(n,d;k)\le \left\lfloor \frac{\gauss{n}{1}{q}\cdot A_q(n-1,d;k-1)}{\gauss{k}{1}{q}}\right\rfloor-\delta-1.\]

Moreover, there exists a $q^r$-divisible multiset of points of cardinality $n$ if and only if there are non-negative integers $a_0,\dots,a_r$ with 
$n=\sum_{i=0}^{r} a_is_{i,r}^q$, where $s_{i,r}^q=q^{r-i}\cdot \frac{q^{i+1}-1}{q-1}$.
\end{Theorem}
This is implemented as \texttt{improved\_johnson}. The iterated version is given by:
\begin{equation}
\label{ie_best_upper_bound}
A_q(n,d;k)\le\left\{ \frac{q^{n}\!-\!1}{q^{k}\!-\!1} \left\{ \frac{q^{n\!-\!1}\!-\!1}{q^{k\!-\!1}\!-\!1} \left\{ \ldots
\left\{ \frac{q^{n'\!+\!1}\!-\!1}{q^{d'\!+\!1}\!-\!1} A_q(n',d;d') \right\}_{d'+1}
\ldots \right\}_{k-2} \right\}_{k-1} \right\}_k,
\end{equation}
where $d'=d/2$, $n'=n-k+d'$, and $\left\{a / \gauss{k}{1}{q}\right\}_k:=b$ with maximal $b\in\mathbb{N}$ permitting 
a representation of $a-b\cdot \gauss{k}{1}{q}$ as non-negative integer combination of the summands $q^{k-1-i}\cdot\frac{q^{i+1}-1}{q-1}$ 
for $0\le i\le k-1$.\footnote{As an example we consider $A_2(9;6;4)\le 
\left\{\gauss{9}{1}{q}A_2(8,6;3)/\gauss{4}{1}{q}\right\}_4=
\left\{\frac{17374}{15}\right\}_4$ using $A_2(8,6;3)=34$. We have $\left\lfloor \frac{17374}{15} \right\rfloor=1158$,  
$17374-1158\cdot 15=4$, $17374-1157\cdot 15=19$, and $17374-1156\cdot 15=34$. Since $4$ and $19$ cannot be written as 
a non-negative linear combination of $8$, $12$, $14$, and $15$, but $34=14+12+8$, we have $A_2(9;6;4)\le 1156$, which improves 
upon the iterative Johnson bound by two. We remark that \cite{kiermaier2017improvement} contains an easy and fast algorithm to check 
the presentability as non-negative integer combination as specified above.}

\subsubsection{MRD bound}
\label{subsec_mrd_bound}

Since the size of the lifted MRD code, see Theorem~\ref{thm_MRD_size}, is quite competitive, it is interesting to 
compare the best known constructions with this very general explicit construction. Even more, lifted MRD codes are 
the basis for more involved constructions, see Subsection~\ref{subsec_echelon_ferrers}. From this point of view it is very 
interesting that an upper bound for the cardinality of constant dimension codes containing the
lifted MRD code (of shape $k \times (n-k)$ and rank distance $d/2$) can be stated:

\begin{Theorem}(see \cite[Theorem 10 and 11]{etzion2013problems})
  Let $\mathcal{C}$ be a constant dimension code with given parameters $q$, $n$, $d$, and $k$ that contains 
  a lifted MRD code. Then:
  \begin{itemize}
    \item if $d=2(k-1)$ and $k \ge 3$ then $\left|\mathcal{C}\right| \le q^{2(n-k)} + A_q(n-k,2(k-2);k-1)$;
    \item if $d=2k$ then $\left|\mathcal{C}\right| \le q^{(n-2k)(k+1)} +
             \gauss{n-2k}{k}{q}\frac{q^n-q^{n-2k}}{q^{2k}-q^k} + A_q(n-2k,2k;2k)$.
  \end{itemize}
\end{Theorem}

\begin{Theorem}[{\cite[Proposition 1]{heinlein2017mrdbound}}]
For $2 \le d/2 \le k \le n-k$ let $C$ be a $(n,\#C,d;k)_q$ CDC that contains an LMRD code.

If $k<d \le 2/3 \cdot n$ we have
\[\#C \le q^{(n-k)(k-d/2+1)} + A_q(n-k,2(d-k);d/2)\text{.}\]
If additionally $d=2k$, $r \equiv n \mod{k}$, $0 \le r <k$, and $\gauss{r}{1}{q}<k$, then the right hand side is equal to $A_q(n,d;k)$ and achievable in all cases.

If $(n, d, k) \in \{ (6+3l,4+2l,3+l), (6l,4l,3l) \mid l \ge 1 \}$, then there is a CDC containing an LMRD with these parameters whose cardinality achieves the bound.

If $k<d$ and $n<3d/2$ we have
\[\#C \le q^{(n-k)(k-d/2+1)} + 1\]
and this cardinality is achieved.

If $d \le k < 3d/2$ we have
\begin{align*}
\#C
&
\le
q^{(n-k)(k-d/2+1)} + A_q(n-k,3d-2k;d)
\\
&
+
\gauss{n-k}{d/2}{q}\gauss{k}{d-1}{q}
q^{(k-d+1)(n-k-d/2)}
/
\gauss{k-d/2}{d/2-1}{q}
\text{.}
\end{align*}
\end{Theorem}

\subsubsection{Bounds for partial spreads}
\label{subsec_upper_bounds_spreads}

Partial spreads attain the maximum possible subspace distance $d=2k$ for constant dimension codes with codewords of dimension $k$. So, it 
does not surprise that good bounds are known for this special case. In this context it makes sense to write $n=tk+r$, where $0\le r<k$.  
The cases $r=0$ and $r=1$ are completely resolved: 

\begin{Theorem}
  \label{thm_spread}(\cite{segre1964teoria}; see also \cite{andre1954nicht}, \cite[p.~29]{dembowski2012finite}, Result 2.1 in \cite{beutelspacher1975partial})
  $\mathbb{F}_q^n$ contains a $k$-spread if and only if $k$ divides $n$, where we assume $1\le k\le n$ and $k,n\in\mathbb{N}$.
\end{Theorem}
The corresponding exact value is implemented as upper bound \texttt{spread}.

\begin{Theorem}
  \label{thm_almost_spread}(\cite{beutelspacher1975partial}; see also \cite{hong1972general} for the special case $q=2$)
  For positive integers $1\le k\le n$ be positive integers with $n\equiv 1\pmod k$
  we have $A_q(n,2k;k)=\frac{q^n-q}{q^k-1}-q+1=q\cdot \frac{q^{n-1}-1}{q^k-1}-q+1=\frac{q^n-q^{k+1}+q^k-1}{q^k-1}$.
\end{Theorem}
The corresponding exact value is implemented as upper bound \texttt{partial\_spread\_2}.

\medskip

Since $\mathbb{F}_q^n$ contains $\gauss{n}{1}{q}=\frac{q^n-1}{q-1}$ points and each $k$-dimensional codeword contains 
$\gauss{k}{1}{q}=\frac{q^k-1}{q-1}$ point, we have:
\begin{Theorem}
\label{spread_bound}
$A_q(n,2k;k) \le \left\lfloor\frac{q^n-1}{q^k-1}\right\rfloor$
\end{Theorem}
This is implemented as \texttt{spread\_bound} and is equivalent to 
Theorem~\ref{thm_johnson_I}. It is tight if and only if $r=0$, where it then matches Theorem~\ref{thm_spread}.    

\begin{Theorem}[\cite{etzion2011error}]
$d=2k \land k \nmid n \Rightarrow A_q(n,d;k) \le \left\lfloor\frac{q^n-1}{q^k-1}\right\rfloor-1 $
\end{Theorem}
This is implemented as \texttt{partial\_spread\_5}. We remark that tighter bounds are known if either $r>1$ or $q>2$, i.e., it is tight 
for $(r,q)=(1,2)$, where it matches Theorem~\ref{thm_almost_spread}. Given the trivial upper bound of Theorem~\ref{spread_bound}, one 
defines $A_q(n,2k;k)=\left\lfloor\frac{q^n-1}{q^k-1}\right\rfloor-\sigma$, where $\sigma$ is called the $\textit{deficiency}$. In these 
terms, we have $\sigma=0$ iff $r=0$ and $\sigma=q-1$ if $r=1$.   

\medskip

For $q=2$ and $k=3$, then requiring $r\in\{0,1,2\}$, the value of $A_2(n,6;3)$ can be determined exactly:
\begin{Theorem} (see \cite{spreadsk3})
  \label{thm_spread_k_3}
  For each integer $m\ge 2$ we have
  \begin{itemize}
    \item[(a)] $A_2(3m,6;3)=\frac{2^{3m}-1}{7}$;
    \item[(b)] $A_2(3m+1,6;3)=\frac{2^{3m+1}-9}{7}$;
    \item[(c)] $A_2(3m+2,6;3)=\frac{2^{3m+2}-18}{7}$.
  \end{itemize}
\end{Theorem}
The corresponding upper bound is implemented as \texttt{partial\_spread\_1}. We remark that it is sufficient to construct 
a code matching $A_2(8,6;3)=34$, which was found by a computer search in \cite{spreadsk3}, to conclude 
$A_2(3m+2,6;3)=\frac{2^{3m+2}-18}{7}$ for all $m\ge 2$, since $\sigma$ is a non-increasing function in $n$, see 
\cite[Lemma~4]{honold2016partial}. The other cases are special instances of $r=0$ and $r=1$.  

For $r=2$ there are the following results:
\begin{Theorem} (Theorem 4.3 in \cite{kurzspreads})
  \label{thm_spread_exact_value_3}
  For each pair of integers $t\ge 1$ and $k\ge 4$ we have $A_2(k(t+1)+2,2k;k)=\frac{2^{k(t+1)+2}-3\cdot 2^{k}-1}{2^k-1}$.
\end{Theorem}
The corresponding upper bound is implemented as \texttt{partial\_spread\_kurz\_q2}.

\begin{Lemma} (Lemma 4.6 in \cite{kurzspreads})
  \label{lemma_spread_upper_bound_3_q}
  For integers $t\ge 1$ and $k\ge 4$ we have $A_3(k(t+1)+2,2k;k) \le \frac{3^{k(t+1)+2}-3^2}{3^k-1}-\frac{3^2+1}{2}$.
\end{Lemma}
The corresponding upper bound is implemented as \texttt{partial\_spread\_kurz\_q3}. We remark that the above two theorems 
are improvements over the following general upper bound from 1979:

\begin{Theorem} (Corollary 8 in \cite{nets_and_spreads})
  \label{thm_partial_spread_4}
  If $n=k(t+1)+r$ with $0<r<k$, then 
\[
    A_q(n,2k;k)\le \sum_{i=0}^{t} q^{ik+r} -\left\lfloor\theta\right\rfloor-1
    =q^r\cdot \frac{q^{k(t+1)}-1}{q^k-1}-\left\lfloor\theta\right\rfloor-1,
\]
  where $2\theta=\sqrt{1+4q^k(q^k-q^r)}-(2q^k-2q^r+1)$.
\end{Theorem}

We remark that this theorem is also restated as Theorem~13 in \cite{etzion2013problems} and as Theorem~44 in \cite{etzionsurvey} 
with the small typo of not rounding down $\theta$ ($\Omega$ in their notation).  
The corresponding upper bound is implemented as \texttt{DrakeFreeman}.

Not too long ago it was shown that the construction of Observation~\ref{obs_multi_component}, i.e., the Echelon-Ferrers construction with 
a skeleton code of disjoint codewords, gives the optimal value if $r\ge 1$ and $k$ is \textit{large enough}, i.e., it is asymptotically 
optimal:   

\begin{Theorem} (see \cite[Theorem~5]{nastase2016maximum})
\label{thm_partial_spread_asymptotic}
For $r=n\pmod k$ and $k>\gauss{r}{1}{q}$ we have:
\[A_q(n,2k;k) = \frac{q^n-q^{k+r}}{q^k-1} +1\]
\end{Theorem}
This is implemented as \texttt{partial\_spread\_NS}. Theorem~\ref{thm_spread_exact_value_3} is just a very special case of it. 
If $k=\gauss{r}{1}{q}$ similar techniques allow to obtain an improved upper bound:
\begin{Theorem}[{\cite[Lemma 10 and Remark 11]{nastase2016maximum}}]
\label{thm_remark11_nastase}
$r=n\pmod k \land k=\gauss{r}{1}{q} < n \land r\ge 2 \Rightarrow A_q(n,2k;k) \le lq^k + \min\{q,\lceil q^r/2\rceil\}$ where 
$l=\frac{q^{n-k}-q^r}{q^k-1}$
\end{Theorem}
This is implemented as \texttt{partial\_spread\_NS\_upper\_bound}.

\medskip

Invoking a result on the existence of so-called \textit{vector space partitions}, see \cite[Theorem 1]{heden2009length}, 
the authors of \cite{nastase2016maximum} obtained the following tightenings:
\begin{Theorem}[{\cite[Theorem 6]{nastase2016maximumII}}]
\label{thm_nastaseII_1}
$r=n\pmod t \land 2 \le r < t \le \gauss{r}{1}{q} \Rightarrow A_q(n,2t;t) \le \frac{q^n-q^{t+r}}{q^t-1} +q^r-(q-1)(t-2)-c_1+c_2$ where $c_1 = 2-t \pmod{q}$ and $ c_2 = \begin{cases} q & q^2 \mid (q-1)(t-2)+c_1 \\ 0 & \text{else} \end{cases} $ such that $-q+1 \le -c_1+c_2 \le q$
\end{Theorem}
This is implemented as \texttt{partial\_spread\_NS\_2\_Theorem6}.

\begin{Theorem}[{\cite[Theorem 7]{nastase2016maximumII}}]
\label{thm_nastaseII_2}
$r=n\pmod t \land 2 \le r < t \le 2^r-1 \Rightarrow A_2(n,2t;t) \le \frac{2^n-2^{t+r}}{2^t-1} +2^r-t+1+c$ where $ c = \begin{cases} 1 & 4 \mid t-1 \\ 0 & \text{else} \end{cases} $
\end{Theorem}
This is implemented as \texttt{partial\_spread\_NS\_2\_Theorem7}.

\medskip

The four bounds mentioned above can possibly be best explained using the concept of divisible codes, see \cite{honold2016partial}. To 
this end we call a point that is not contained in any $k$-dimensional codeword of a partial spread a \textit{hole}. Taking the set of 
holes as columns of a generator matrix, 
we obtain a projective linear code over $\mathbb{F}_q$ of dimension $k$ and the number of holes as length. It turns out that the Hamming 
weights of all codewords are divisible by $q^{k-1}$, see \cite[Theorem 8]{honold2016partial}. Those codes have to satisfy the famous 
\textit{MacWilliams Identities}, see \cite{macwilliams63},
\begin{equation}
  \label{mac_williams_identies}
  \sum_{j=0}^{n-i} {{n-j}\choose i} A_j=q^{k-i}\cdot \sum_{j=0}^i
  {{n-j}\choose{n-i}}A_j^\perp\quad\text{for }0\le i\le n, 
\end{equation} 
where $A_j$ denotes the number of codewords of Hamming weight $j$ and $A_j^\perp$ denotes the number of codewords of weight $j$ of the dual code. 
We have $A_0=A_0^\perp=1$ and $A_1^\perp=0$. Projectivity of the code is equivalent to $A_2^\perp=0$ and the divisibility conditions says that  
$A_j=0$ for all indices $j$ that are not divisible by $q^{k-1}$. Moreover, the \textit{residual codes} are $q^{k-2}$-divisible, which can be 
applied recursively. The first two MacWilliams Identities can be used to exclude the existence of quite some lengths of projective linear
$q^r$-divisible codes. Translated back to partial spreads, this gives:
\begin{Theorem}[\cite{kurz2016upper,kurz2017packing}]
\label{ps_two_mac_williams}
$r \ge 1 \land k \ge 2 \land z,u \ge 0 \land t = \gauss{r}{1}{q} +1 -z+u > r \Rightarrow A_q(n,2t;t) \le lq^t+1+z(q-1)$ where $l=\frac{q^{n-t}-q^r}{q^t-1}$ and $n=kt+r$
\end{Theorem}
This is implemented as \texttt{partial\_spread\_kurz16\_28} and contains Theorem~\ref{thm_partial_spread_asymptotic} and 
Theorem~\ref{thm_remark11_nastase} as a special case. If the infeasibility of the first three MacWilliams Identities is used, one obtains:
\begin{Theorem}[{\cite[Theorem 2.10]{kurz2017packing},\cite[Theorem 10]{honold2016partial}}]
\label{ps_three_mac_williams}
$r \ge 1 \land t \ge 2 \land y \ge \max\{r,2\} \land z \ge 0 \land r,t,y,z \in \mathbb{Z} \land u=q^y \land y \le k \land k= \gauss{r}{1}{q} +1 -z >r \land n=kt+r \land l=\frac{q^{n-k}-q^r}{q^k-1}$ $\Rightarrow A_q(n,2k;k) \le lq^k + \lceil u-1/2 - 1/2\sqrt{1+4u(u-(z+y-1)(q-1)-1)}\rceil$ Note that the description contains the value of $y$ in brackets
\end{Theorem}
This is implemented as \texttt{partial\_spread\_HKK16\_T10}. Setting $y=k$ in Theorem~\ref{ps_three_mac_williams} gives 
Theorem~\ref{thm_partial_spread_4}, i.e., the \textit{classical} result of Drake and Freeman. We remark that the combination of 
Theorem~\ref{ps_two_mac_williams} and Theorem~\ref{ps_three_mac_williams} is at least as tight as the combination of   
Theorem~\ref{thm_partial_spread_asymptotic}, Theorem~\ref{thm_remark11_nastase}, Theorem~\ref{thm_nastaseII_1} and Theorem~\ref{thm_nastaseII_2} and 
in several cases the first mentioned two theorems are strictly tighter. This statement was numerically verified for all 
$2\le q\le 9$, $1\le n,k\le 100$ in \cite{honold2016partial}. There is also a conceptual reason: The result of Heden on the existence 
of vector space partitions, see \cite[Theorem 1]{heden2009length}, can be improved by using the implications of the first three MacWilliams Identities 
for divisible codes, see \cite[Theorem 12]{honold2016partial}, which classifies the possible length $n$ of $q^r$-divisible codes for all $n\le rq^{r+1}$. 
For larger $n$ some partial numerical results are obtained in \cite{heinlein2017projective}. A further, more direct, improvement of Heden's result can 
be found in \cite{kurz2017heden}.

Excluding codes by showing that the Equation~(\ref{mac_williams_identies}) has no non-negative real solution is known as the 
linear programming method, which generally works for association schemes, see \cite{delsarte1973algebraic}. For Theorem~\ref{ps_two_mac_williams} and 
Theorem~\ref{ps_three_mac_williams} only a first few equations are taken into account and an analytical solution was obtained. For the first 
four equations of (\ref{mac_williams_identies}) the following analytical criterion was stated in \cite{honold2016partial}:  

\begin{Lemma}
  \label{lemma_implication_fourth_mac_williams}
  Let $\mathcal{C}$ be $\Delta$-divisible over $\mathbb{F}_q$ of cardinality
  $n>0$ and $t\in\mathbb{Z}$. Then
  $ \sum_{i\ge 1} \Delta^2(i-t)(i-t-1)\cdot (g_1\cdot i+g_0)\cdot
  A_{i\Delta}\,\,+qhx = n(q-1)(n-t\Delta)(n-(t+1)\Delta)g_2 $, where
  $g_1=\Delta qh$, $g_0=-n(q-1)g_2$,
  $g_2=h-\left(2\Delta qt+\Delta q-2nq+2n+q-2\right)$ and
  $ h= \Delta^2q^2t^2+\Delta^2q^2t-2\Delta nq^2t-\Delta nq^2+2\Delta
  nqt+n^2q^2+\Delta nq-2n^2q+n^2+nq-n $.
\end{Lemma}

\begin{Corollary}
  \label{cor_implication_fourth_mac_williams}
  If there exists $t\in\mathbb{Z}$, 
  using the notation of
  Lemma~\ref{lemma_implication_fourth_mac_williams}, with
  $n/\Delta\notin [t,t\!+\!1]$, $h\ge 0$, and $g_2<0$, then 
  there
  is 
  no $\Delta$-divisible
  set over $\mathbb{F}_q$ of cardinality $n$.
\end{Corollary}

Numerically evaluating this criterion or numerically solving the corresponding linear programs, 
and taking into account that the deficiency $\sigma$ is non-increasing in $n$, gives:
\begin{Theorem}[\cite{kurz2016upper,kurz2017packing}]
\label{ps_four_mac_williams}
$A_2(4k+3,8;4)\le 2^4l+4$, where $l=\frac{2^{4k-1}-2^3}{2^4-1}$ and $k \ge 2$, \\
$A_2(6k+4,12;6)\le 2^6l+8$, where $l=\frac{2^{6k-2}-2^4}{2^6-1}$ and $k \ge 2$, \\
$A_2(6k+5,12;6)\le 2^6l+18$, where $l=\frac{2^{6k-1}-2^5}{2^6-1}$ and $k \ge 2$, \\
$A_3(4k+3,8;4)\le 3^4l+14$, where $l=\frac{3^{4k-1}-3^3}{3^4-1}$ and $k \ge 2$, \\
$A_3(5k+3,10;5)\le 3^5l+13$, where $l=\frac{3^{5k-2}-3^5}{3^3-1}$ and $k \ge 2$, \\
$A_3(5k+4,10;5)\le 3^5l+44$, where $l=\frac{3^{5k-1}-3^4}{3^5-1}$ and $k \ge 2$, \\
$A_3(6k+4,12;6)\le 3^6l+41$, where $l=\frac{3^{6k-2}-3^4}{3^6-1}$ and $k \ge 2$, \\
$A_3(6k+5,12;6)\le 3^6l+133$, where $l=\frac{3^{6k-1}-3^5}{3^6-1}$ and $k \ge 2$, \\
$A_3(7k+4,14;7)\le 3^7l+40$, where $l=\frac{3^{7k-3}-3^4}{3^7-1}$ and $k \ge 2$, \\
$A_4(5k+3,10;5)\le 4^5l+32$, where $l=\frac{4^{5k-2}-4^3}{4^5-1}$ and $k \ge 2$, \\
$A_4(6k+3,12;6)\le 4^6l+30$, where $l=\frac{4^{6k-3}-4^3}{4^6-1}$ and $k \ge 2$, \\
$A_4(6k+5,12;6)\le 4^6l+548$, where $l=\frac{4^{6k-1}-4^5}{4^6-1}$ and $k \ge 2$, \\
$A_4(7k+4,14;7)\le 4^7l+128$, where $l=\frac{4^{7k-3}-4^4}{4^7-1}$ and $k \ge 2$, \\
$A_5(5k+2,10;5)\le 5^5l+7$, where $l=\frac{5^{5k-3}-5^2}{5^5-1}$ and $k \ge 2$, \\
$A_5(5k+4,10;5)\le 5^5l+329$, where $l=\frac{5^{5k-1}-5^4}{5^5-1}$ and $k \ge 2$, \\
$A_7(5k+4,10;5)\le 7^5l+1246$, where $l=\frac{7^{5k-1}-7^2}{7^5-1}$ and $k \ge 2$, \\
$A_8(4k+3,8;4)\le 8^4l+264$, where $l=\frac{8^{4k-1}-8^3}{8^4-1}$ and $k \ge 2$, \\
$A_8(5k+2,10;5)\le 8^5l+25$, where $l=\frac{8^{5k-3}-8^2}{8^5-1}$ and $k \ge 2$, \\
$A_8(6k+2,12;6)\le 8^6l+21$, where $l=\frac{8^{6k-4}-8^2}{8^6-1}$ and $k \ge 2$, \\
$A_9(3k+2,6;3)\le 9^3l+41$, where $l=\frac{9^{3k-1}-9^2}{9^3-1}$ and $k \ge 2$, and \\
$A_9(5k+3,10;5)\le 9^5l+365$, where $l=\frac{9^{5k-2}-9^3}{9^5-1}$ and $k \ge 2$
\end{Theorem}
This is implemented as \texttt{partial\_spread\_kurz16\_additional}. We remark that we are not aware 
of a set of parameters, where considering more than four equations from (\ref{mac_williams_identies}) yields 
an improvement.

\subsubsection{Further upper bounds}
\label{subsec_further_upper}

\begin{Theorem}[{\cite[Upper bound of Theorem 3.3(i)]{MR3543542}}] 
\(A_q(6,4)\le(q^3+1)^2\) for all \(q\ge 3\).
\end{Theorem}

This is implemented as \texttt{HKK\_theorem\_3\_3\_i\_upper\_bound\_A\_q(6,4)}.

\begin{Theorem}[{\cite[Upper bound of Lemma 2.4]{MR3543542}}] 
\(1\le d/2\leq k\le\lfloor n/2\rfloor, d\equiv 0\bmod 2 \Rightarrow A_q(n,d;k)>q\cdot A_q(n,d;k-1)\).
\end{Theorem}

This is implemented as \texttt{HKK\_theorem\_3\_3\_i\_upper\_bound\_A\_q(6,4)}.

\begin{Theorem}
(see \cite[Theorem~3]{ahlswede2009error})
\label{thm_ahlswede_aydinian}
For $0\le t< r\le k$, $k-t\le m\le v$, and $t\le v-m$ we have:
\[A_q(n,2r;k)\le \frac{\gauss{n}{k}{q} A_q(m,2r-2t;k-t)}{\sum_{i=0}^t q^{i(m+i-k)}\gauss{m}{k-i}{q}\gauss{n-m}{i}{q}}\]
\end{Theorem}
Note that the description of the application of the constraint contains $t$, $m$, and an optional $o$, indicating the application on the parameters of the orthogonal 
code, in brackets. This bound is implemented as \texttt{Ahlswede\_Aydinian}. We remark that there are typos in the formulation in 
\cite{ahlswede2009error,khaleghi2009subspace}. The corrected stated version can also be found in \cite{heinlein2017asymptotic}. 
The authors of \cite{ahlswede2009error} have observed that Theorem~\ref{thm_ahlswede_aydinian} contains Theorem~\ref{thm_johnson_cdc}, 
i.e., the Johnson bound, as a special case. In \cite{heinlein2017asymptotic} it was numerically checked that Theorem~\ref{thm_ahlswede_aydinian} 
does not give strictly tighter bounds than Theorem~\ref{thm_johnson_cdc} for all $2 \le q \le 9$, $4 \le n \le 100$, and $4 \le d \le 2k \le n$.  

\medskip

The Delsarte linear programming bound for the $q$-Johnson scheme, which is an association scheme, was obtained in \cite{delsarte1978hahn}. However, 
numerical computations indicate that it is not better than the Anticode bound, see \cite{MR3063504}. In \cite{zhang2011linear} it was 
shown that the Anticode bound is implied by the Delsarte linear programming bound. In \cite{MR3063504} it was shown that a semidefinite programming 
formulation, that is equivalent to the Delsarte linear programming bound, implies the Anticode bound of Theorem~\ref{theo:anticode}, the sphere-packing 
of Theorem~\ref{thm:sphere_packing}, the \textit{weak} Johnson bound of Theorem~\ref{thm_johnson_I}, and the Johnson bound of  
Theorem~\ref{thm_johnson_cdc} (without rounding). This makes perfectly sense, since Theorem~\ref{thm:sphere_packing} and Theorem~\ref{thm_johnson_I} 
are implied by Theorem~\ref{theo:anticode} and the iteration of Theorem~\ref{thm_johnson_cdc} without rounding gives exactly Theorem~\ref{theo:anticode}. 
Using \texttt{Maple} and exact arithmetic, we have checked that for all $2\le q\le 9$, $4\le n\le 19$, $2\le k\le n/2$, $4\le d\le 2k$ the optimal value of the Delsarte linear programming 
bound is indeed the Anticode bound. Given the result from \cite{zhang2011linear} it remains to construct a feasible solution of the 
Delsarte linear programming formulation whose target value equals the Anticode bound. Such a feasible solution can also be constructed recursively. 
To this end, let $x_0$, \dots, $x_{k-1}$ denote a primal solution for the parameters of $A_q(n-1, d, k-1)$, then $z_0$, \dots, $z_k$ 
is a feasible solution for the parameters of $A_q(n, d,k)$ setting $z_i=x_i \cdot  \gauss{k}{1}{q} \ \gauss{k-i}{1}{q}$ for all
$0 \le i \le k-1$ and $z_k= \gauss{n}{k}{q}/\gauss{n-k+d/2-1}{d/2-1}{q} - z_0 -  \dots -z_{k-1}$. For the mentioned parameter space this conjectured 
primal solution is feasible with the Anticode bound as target value. Due to the property of the symmetry group of 
$(\mathbb{F}_q^n,d_S)$, i.e., two-point homogeneous, the symmetry reduced version of the semidefinite programming formulation of the maximum clique 
problem formulation collapses the Delsarte linear programming bound for the $q$-Johnson scheme.

\medskip

Another rather general technique to obtain upper bounds for the maximum cliques size of a graph is to use $p$-ranks, see e.g.\ \cite[Lemma~1.3]{MR3692294}.
\begin{Lemma}
  Let $G$ be a graph with adjacency matrix $A$ and $Y$ be a clique of $G$, then
  $$
    |Y|\le\left\{
    \begin{array}{rcl}
      \operatorname{rank}_p(A)+1 && \text{if $p$ divides $|Y|-1$,}\\
      \operatorname{rank}_p(A)   & & \text{otherwise.}
    \end{array}
    \right.
  $$
\end{Lemma}
This is implemented as \texttt{prank}. Some numerical experiments suggest that the resulting upper bounds are rather weak, e.g.,
$A_2(4,4;2) \le 5$,
$A_2(5,4;2) \le 19$,
$A_2(6,4;2) \le 49$,
$A_2(6,4;3) \le 223$, and
$A_2(6,6;3) \le 19$.

\medskip

We close this section by upper bounds obtained from tailored integer linear programming 
computations. The five optimal isomorphism types for $A_2(6,4;3)=77$ have been determined in 
\cite{hkk77}. The upper bound $A_2(8,6;4) \le 272$ was obtained in \cite{heinlein2017new} and 
is implemented as \texttt{special\_case\_2\_8\_6\_4}. A little later the exact value 
$A_2(8,6;4)=257$ and its two optimal isomorphism types were determined, see \cite{heinlein2017classifying}. 
The assumptions can even be weakened and still yield an upper bound strictly less than $289$, which 
follows from the Johnson bound, see \cite{heinlein2019generalized_vsp}. 

\subsubsection{Exact Bounds}
\begin{Theorem}[{\cite[Theorem 3.1(ii)]{MR3543542}}] 
If \(n=2k\) is even then \(A_q(n,n;k)=q^k+1\).
\end{Theorem}

This is implemented as \texttt{HKK\_theorem\_3\_1\_ii\_cdc}.

\begin{Theorem}[{\cite[Theorem 3.2(ii)]{MR3543542}}] 
If \(n=2k+1\ge 5\) is odd then \(A_q(n,n-1;k)=q^{k+1}+1\).
\end{Theorem}

This is implemented as \texttt{HKK\_theorem\_3\_2\_ii\_cdc}.

\section{Bounds for MDCs}
\label{sec_mdc}

For mixed dimension subspace codes the choice between the subspace distance $d_S$ and the injection distance 
$d_I$ really makes a difference. Here we consider the subspace distance only. Subsection~\ref{subsec_mdc_lower} 
is devoted to constructions and upper bounds are presented in Subsection~\ref{subsec_mdc_upper}. In general, mixed 
dimension subspace codes have obtained much less attention than constant dimension codes. Obtaining bounds seems to be 
more challenging. For surveys we refer e.g.\ to \cite{MR3063504,MR3543542,khaleghi2009subspace}. 

\subsection{Lower bounds and constructions}
\label{subsec_mdc_lower}
Of course the empty set is a mixed dimension code for any dimension and subspace distance. 
$A_q(n,d)\ge 0$ is implemented as \texttt{trivial\_2}. If $d\le 2n$ and $n\ge 1$, then 
$\{\langle 0 \rangle , \mathbb{F}_q^n\}$ is a mixed dimension codes, so that $A_q(n,d)\ge 2$. 
This is implemented as \texttt{trivial\_4}. We structure the following lower bounds and constructions into  
nonrecursive, see Subsection~\ref{subsubsec_nonrecursive_lower_mdc}, and recursive lower bounds, 
see Subsection~\ref{subsubsec_recursive_lower_mdc}. Some constraints leading to exact values are 
also collected in Subsection~\ref{subsubsection_exact_value_mdc}. A survey on lower bounds and constructions 
can be found in \cite{horlemann2018constructions}. 

\subsubsection{Nonrecursive lower bounds}
\label{subsubsec_nonrecursive_lower_mdc}

The Echelon-Ferrers construction also works for mixed dimension codes, see e.g.\ \cite{etzion2009error,gorla2014subspace}. 
Also the stated ILP formulation directly transfers, which is implemented as \texttt{echelon\_ferrers}. A more sophisticated 
search for the optimal construction within this setting is implemented as \texttt{ef\_computation}, 
cf.\ Subsection~\ref{subsec_echelon_ferrers}.

Similar to the sphere covering bound for constant dimension codes in Theorem~\ref{thm_cdc_sphere_covering}, there 
exists a version for mixed dimension codes.:
\begin{Theorem}[{\cite[Theorem 9]{etzion2011error}}]
\[
A_q(n,d) \ge
\frac{
\sum_{k=0}^n \sum_{j=0}^n \gauss{n}{k}{q}\gauss{n}{j}{q}
}{
\sum_{k=0}^n \sum_{j=0}^{d-1} \sum_{i=0}^{j}
\gauss{n}{k}{q}\gauss{k}{i}{q}\gauss{n-k}{j-i}{q} q^{i(j-i)}
}
\]
\end{Theorem}
This is implemented as \texttt{gilbert\_varshamov}.

\begin{Theorem}[{\cite[Theorem 3.3(ii)]{MR3543542}}]
$ A_q(n,n-2) \ge 2q^{k+1}+1 $ for $ n=2k+1 \ge 5 $
\end{Theorem}
This is implemented as \texttt{nodd\_deqnm2\_l}.

\subsubsection{Recursive lower bounds}
\label{subsubsec_recursive_lower_mdc}

\begin{Theorem}[\cite{etzion2009error}]
$ \left\lceil \max_{k=0}^n \frac{q^{n+1-k}+q^k-2}{q^{n+1}-1} \cdot A_q(n+1,d+1;k) \right\rceil \le A_q(n,d)$
\end{Theorem}
This is implemented as \texttt{cdc\_average\_argument}.

\begin{Theorem}
$\max_{k=0}^n A_q(n,d;k) \le A_q(n,d)$
\end{Theorem}
This is implemented as \texttt{cdc\_lower\_bound}.

\begin{Theorem}[{\cite[Lower bound of Theorem 3.3(i)]{MR3543542}}]
If $n=2k \ge 8$ even then $A_q(n,n-2)\ge A_q(n,n-2,k)$.
\end{Theorem}
This is implemented as \texttt{HKK\_theorem\_3\_3\_i\_lower\_bound}.

\begin{Theorem}[{\cite[Lower bound of Theorem 2.5]{MR3543542}}] \label{theo:improved_cdc_lower_bound}
$ \sum_{k=0 \land k \equiv \lfloor n / 2 \rfloor \pmod{d}}^{n} A_q(n,2\lceil d/2 \rceil;k) \le A_q(n,d) \le 2 + \sum_{k=\lceil d/2 \rceil}^{n- \lceil d/2 \rceil} A_q(n,2\lceil d/2 \rceil;k) $
\end{Theorem}
This is implemented as \texttt{improved\_cdc\_lower\_bound}.

\begin{Theorem}[\cite{MR3543542}]
The bound is $ A_q(n,d) \ge \max\{ \sum_{k \in K} A_q(n,d;k) \mid K \subseteq \{0,\ldots,n\} : |k_1-k_2| \ge d \;\forall k_1 \ne k_2 \in K \}$. This is computed using dynamic programming and the function $ L(N) := \max\{\sum_{k \in K} A_q(n,d;k) \mid K \subseteq \{0,\ldots,N\} : |k_1-k_2| \ge d \;\forall k_1 \ne k_2 \in K\} = \max\{ L(N-1), L(N-d) + A_q(n,d;N) \} $ for all $N = 0,\ldots,n$.
\end{Theorem}
This is implemented as \texttt{layer\_construction}.

\begin{Theorem}[{\cite[Lower bound of Theorem 3.3(i)]{MR3543542}}]
\(q^6+2q^2+2q+1\le A_q(6,4)\) for all \(q\ge 3\).
\end{Theorem}
This is implemented as \texttt{HKK\_theorem\_3\_3\_i\_lower\_bound\_A\_q(6,4)}.

\subsection{Implemented upper bounds}
\label{subsec_mdc_upper}
Quoting \cite{MR3063504}, bounds for mixed dimension codes are much harder to obtain than for 
constant dimension codes, since, for example, the size of balls in this space depends not only 
on their radius, but also on the dimension of their center. We structure the upper bound into 
nonrecursive, see Subsection~\ref{subsubsec_nonrecursive_upper_mdc}, and recursive bound, see 
Subsection~\ref{subsubsec_recursive_upper_mdc}. Some constraints leading to exact values are 
also collected in Subsection~\ref{subsubsection_exact_value_mdc}.

\subsubsection{Nonrecursive upper bounds}
\label{subsubsec_nonrecursive_upper_mdc}

\begin{Theorem}
If $d=n$ then the whole vector space is the direct sum of each pair of codewords. If a code had three codewords, 
then $2k=n$ which is impossible for $n$ odd.
\end{Theorem}
This is implemented as \texttt{nodd\_deqn}.

\begin{Theorem}[{\cite[Upper bound of Theorem 3.3(ii)]{MR3543542}}]
$ A_q(n,n-2) \le 2q^{k+1}+2 $ for $ n=2k+1 \ge 5 $.
\end{Theorem}
This is implemented as \texttt{nodd\_deqnm2\_u}.

\medskip

The Johnson bound for constant dimension codes can be modified for mixed dimension subspace codes, see \cite{ubt_eref52236}. 
\begin{Theorem}[{\cite[Lemma 2]{ubt_eref52236}}]
$$
  A_2(10,5)\le 48104.
$$
\end{Theorem}
This bound is implemented as \texttt{johnson\_MDC\_Lemma\_2}.

\begin{Theorem}[{\cite[Lemma 3]{ubt_eref52236}}]
$$
  A_3(9,5)\le 123048.
$$
\end{Theorem}
This bound is implemented as \texttt{johnson\_MDC\_Lemma\_3}.

\begin{Theorem}[{\cite[Lemma 5]{ubt_eref52236}}]
$$
  A_q(7,3)\le 2\left(q^8+q^6+2q^5+2q^3+q^2-q+2\right).
$$
\end{Theorem}
This bound is implemented as \texttt{johnson\_MDC\_Lemma\_5}.

\begin{Theorem}[{\cite[Proposition 5]{ubt_eref52236}}]
We have $A_2(8,3)\le 9260$ and $A_q(8,3)\le q^{12}+3q^{10}+q^9+3q^8+3q^7+3q^6+5q^5+3q^4+q^3+4q^2+2q-1$ for $q\ge3$.
\end{Theorem}
This bound is implemented as \texttt{johnson\_MDC\_Proposition\_5}.

\medskip

For the mixed dimension case the acting symmetry group is not $2$-point homogeneous, so that the semidefinite programming formulation 
of the maximum clique problem after symmetrization does not collapse to a linear program. Numerical evaluations of this SLP are given by:
\begin{Theorem}[\cite{MR3063504}]
\label{thm_slp}
$A_2(4,3) \le 6$,
$A_2(5,3) \le 20$,
$A_2(6,3) \le 124$,
$A_2(7,3) \le 776$,
$A_2(7,5) \le 35$,
$A_2(8,3) \le 9268$,
$A_2(8,5) \le 360$,
$A_2(9,3) \le 107419$,
$A_2(9,5) \le 2485$,
$A_2(10,3) \le 2532929$,
$A_2(10,5) \le 49394$,
$A_2(10,7) \le 1223$,
$A_2(11,5) \le 660285$,
$A_2(11,7) \le 8990$,
$A_2(12,7) \le 323374$,
$A_2(12,9) \le 4487$,
$A_2(13,7) \le 4691980$,
$A_2(13,9) \le 34306$,
$A_2(14,9) \le 2334086$,
$A_2(14,11) \le 17159$,
$A_2(15,11) \le 134095$, and
$A_2(16,13) \le 67079$.
\end{Theorem}
This is implemented as \texttt{semidefinite\_programming}. See also \cite{heinlein2018new}.

\begin{Theorem}{\cite[Theorem 4.10]{MR3543542}}
$A_2(6,3) \le 118$ and $A_2(7,4) \le 407$
\end{Theorem}
This is implemented as \texttt{special\_cases\_upper\_notderived}.

\begin{Theorem}
A subspace code is a subset of the subspaces of $\mathbb{F}_q^n$, i.e., $A_q(n,d) \le \sum_{k=0}^n \gauss{n}{k}{q}$.
\end{Theorem}
This is implemented as \texttt{trivial\_3}.

\subsubsection{Recursive upper bounds}
\label{subsubsec_recursive_upper_mdc}

\begin{Theorem}
$A_q(n,d) \le \sum_{k=0}^n A_q(n,d;k)$
\end{Theorem}
This is implemented as \texttt{cdc\_upper\_bound}.

\medskip

The following approach generalizes the sphere-packing bound for constant dimension codes 
facing the fact that the spheres have different sizes. To that end let $B(V,e)$ denote the ball with center $V$ and radius $e$. 
Those balls around codewords are pairwise disjoint.

\begin{Theorem}{\cite[Theorem 10]{etzion2011error}}
\label{thm_ilp_mdc_ev}
Denoting the number of $k$-dimensional subspaces contained in $B(V,e)$ with $\dim(V)=i$ by 
$c(i,k,e)$, we have
\begin{equation*}
  c(i,k,e)=\sum_{j=\left\lceil\frac{i+k-e}{2}\right\rceil}^{\min\{k,i\}} \gauss{i}{j}{q}\gauss{n-i}{k-j}{q} q^{(i-j)(k-j)}.
\end{equation*} 
Thus, $A_q(n,2e+1)$ is at most as large as the target value of:
\begin{align}
\max \sum_{i=0}^n &a_i \label{ILP_EtzionVardy}\\
\text{subject to}\,\, &a_i \le A_q(n,2e+2;i)  && \forall 0\le i\le n \nonumber\\
\sum_{i=0}^n c(i,k,e)\cdot &a_i \le \gauss{n}{k}{q} && \forall 0\le k\le n  \nonumber\\
&a_i \in \mathbb{N} &&\forall 0\le i\le n\nonumber
\end{align}
\end{Theorem}
This is implemented as \texttt{Etzion\_Vardy\_ilp}.

\medskip

\cite[Theorem~10]{khaleghi2009subspace} refers to another LP upper bound by Ahlswede and Aydinian, see \cite{ahlswede2009error}.

\begin{Theorem}[{\cite[Theorem 5 and Theorem 6]{ahlswede2009error},\cite[Theorem 10]{khaleghi2009subspace}}]
For integers \(1\le t \le {n\over 2} \), let
\begin{equation}
 f(n, t, q)=\max(f_{0}+f_{1}+\ldots+f_{n})
\end{equation}
subject to linear constraints:
\begin{equation}
f_{0}, f_{1}, \dots, f_{n}
\end{equation}
are nonnegative integers with
\begin{align}
&f_{0}=f_{n}=1, \, f_{k}=f_{n-k}=0  &\text{ for } k=1, \dots, t \\
&f_{k}\leq A_{q}(n, 2t+2, k) &\text{ for } k=0, \dots, n\\ 
&f_{k}+{1\over t+1}\sum_{i=1}^{t}(t+1-i)(f_{k-i}\gauss{n-k+i}{n-k}{q}+f_{k+i}\gauss{k +i}{k}{q})  \le \gauss{n}{k}{q} &\text{ for } k=0, \dots, n\\
&f_{-j}=f_{n+j} =0 &\text{ for }  i=1, \dots, t.
\end{align}
Then
\begin{equation}
A_{q}(n, 2t+1)\le f(n, t, q).
\end{equation}
\end{Theorem}

This is implemented as \texttt{Ahlswede\_Aydinian\_ilp}.

\begin{Theorem}[{\cite[Upper bound of Theorem 2.5]{MR3543542}}] \label{theo:improved_cdc_upper_bound}
$ \sum_{k=0 \land k \equiv \lfloor n / 2 \rfloor \pmod{d}}^{v} A_q(n,2\lceil d/2 \rceil;k) \le A_q(n,d) \le 2 + \sum_{k=\lceil d/2 \rceil}^{n- \lceil d/2 \rceil} A_q(n,2\lceil d/2 \rceil;k) $
\end{Theorem}
This is implemented as \texttt{improved\_cdc\_upper\_bound}.

\begin{Theorem}[{\cite[Upper bound of Theorem 3.3(i)]{MR3543542}}]
If $n=2k \ge 8$ even then $A_q(n,n-2)\le A_q(n,n-2,k)$.
\end{Theorem}
This is implemented as \texttt{HKK\_theorem\_3\_3\_i\_upper\_bound}.

\begin{Theorem}
$A_q(n,d) \le A_q(n,d-1)$
\end{Theorem}
This is implemented as \texttt{relax\_d}. This innocent and trivial looking inequality produces the tightest known upper bound 
in our database since e.g.\ Theorem~\ref{thm_slp} and Theorem~\ref{thm_ilp_mdc_ev} are not evaluated for all parameters.

\begin{Theorem}[{\cite[Upper bound of Theorem 3.3(i)]{MR3543542}}]
\(A_q(6,4)\le(q^3+1)^2\) for all \(q\ge 3\).
\end{Theorem}
This is implemented as \texttt{HKK\_theorem\_3\_3\_i\_upper\_bound\_A\_q(6,4)}.

\subsection{Further constraints which determine an exact value}
\label{subsubsection_exact_value_mdc}

\begin{Theorem}[{\cite[Theorem 3.4]{MR3543542}}]
If $v=2k$ is even then \(A_q(n,2)=\sum_{0\leq i\leq n\land i\equiv k\bmod 2}\binom{n}{i}_{q}\).
If \(n=2k+1\) is odd then \(A_q(n,2)=\sum_{0\leq i\leq n\land i\equiv0\bmod 2}\binom{n}{i}_{q}=\sum_{0\leq i\leq n\land i\equiv 1\bmod 2}\binom{n}{i}_{q}\).
\end{Theorem}
This is implemented as \texttt{d2}.

\begin{Theorem}[{\cite[Theorem 3.1.ii]{MR3543542}}]
$ A_q(n,n) = q^k+1 $ for $ n=2k $
\end{Theorem}
This is implemented as \texttt{neqdeven}.

\begin{Theorem}[{\cite[Theorem 3.2.i]{MR3543542}}]
$ A_q(n,n-1) = q^k+1 $ for $ n=2k \ge 4 $
\end{Theorem}
This is implemented as \texttt{neven\_deqnm1}.

\begin{Theorem}[{\cite[Theorem 3.2.ii]{MR3543542}}]
$ A_q(n,n-1) = q^{k+1}+1 $ for $ n=2k+1 \ge 5 $
\end{Theorem}
This is implemented as \texttt{nodd\_deqnm1}.

\begin{Theorem}[{\cite[Theorem 3.3(ii)]{MR3543542}}; see also \cite{ghatak2017optimal} and footnote~44 in \cite{MR3543542} referring to 
independent (still unpublished) work of Cossidente, Pavese and Storme]
$A_q(5,3)= 2q^3+2$
\end{Theorem}
This is implemented as \texttt{n5\_d3\_CPS}.

\begin{Theorem}[{\cite[Theorem 3.3(ii)]{MR3543542}}]
$ A_q(5,3) = 2q^3+2 $ for all $q$ and $ A_2(7,5)=34 $
\end{Theorem}
This is implemented as \texttt{nodd\_deqnm2\_e}. We remark that the $20$ isomorphism types of
all latter optimal codes have been classified in \cite{honold2016classification}.

\begin{Theorem}
If the distance is 0 or 1 then the optimal subspace code consists of all subspaces of $\mathbb{F}_q^n$, i.e., $d \le 1 \Rightarrow A_q(n,d) = \sum_{k=0}^n \gauss{n}{k}{q}$.
\end{Theorem}
This is implemented as \texttt{trivial\_dle1}.

\begin{Theorem}[{\cite[Theorem 3.3(i)]{MR3543542}}]
If \(n\) is odd then \(A_q(n,n)=2\).
\end{Theorem}
This is implemented as \texttt{HKK\_theorem\_3\_3\_i}.

\section{Conclusion}
\label{sec_conclusion}

\noindent
The collection of the known results on lower and upper bounds for subspace codes is an ongoing project. So far we have 
merely implemented the tip of the iceberg of the available knowledge. Even for upper bounds for constant dimension codes, 
which is the most advanced part of our summary, there are several pieces of work left. We still hope that the emerging on-line data base 
and the accompanying user's guide is already valuable for researchers in the field at this stage. One yardstick for our knowledge 
is the fraction between the best known lower bound and the best known upper bound for constant dimension codes. To be able to state 
some parametrical results, we compare the size of the lifted MRD code with the Singleton or the Anticode bound as done 
in \cite{heinlein2017asymptotic}. To this end we utilize the so called $q$-Pochhammer symbol
$(a;q)_n:=\prod_{i=0}^{n-1} \left(1-aq^i\right)$ and its specialization $(1/q;1/q)_n=\prod_{i=1}^{n}\left(1-1/q^i\right)$. 
\begin{Proposition}{\cite[Proposition 7]{heinlein2017asymptotic}}
  For $k \le n-k$ the ratio of the size of an LMRD code divided by the size of the Singleton bound converges
  for $n \rightarrow \infty$ monotonically decreasing to $(1/q;1/q)_{k-d/2+1}\ge (1/2;1/2)_\infty > 0.288788$.
\end{Proposition}
\begin{Proposition}{\cite[Proposition 8]{heinlein2017asymptotic}}
  \label{prop_ratio_lmrd_anticode}
  For $k \le n-k$ the ratio of the size of an LMRD code divided by the size of the Anticode bound converges
  for $n \rightarrow \infty$ monotonically decreasing to $\frac{(1/q;1/q)_{k}}{(1/q;1/q)_{d/2-1}}\ge \frac{q}{q-1}\cdot (1/q;1/q)_{k}
  \ge 2\cdot (1/2;1/2)_\infty > 0.577576$.
\end{Proposition}
The largest gap of this estimate is attained for $d=4$ and $k=\left\lfloor n/2\right\rfloor$. We remark that for this special case 
none of the mentioned upper bounds yields an asymptotic improvement over the Anticode bound and none of the described constructions 
yields an asymptotic improvement over the LMRD code construction. If $k$ does not vary with $n$ (or does increase very slowly), then 
the Anticode bound can be asymptotically be attained by an optimal code, see \cite[Theorem~4.1]{MR829351} and also \cite{blackburn2012asymptotic}.

For mixed dimension codes comparatively little is known and more research is sorely needed. If you want to support us in our task, please 
let us know any known constructions, bounds or papers that we have missed so far via  \href{mailto:daniel.heinlein@aalto.fi}{daniel.heinlein@aalto.fi} 
or the \textit{Contribute}-button in the upper right corner of the webpage \url{subspacecodes.uni-bayreuth.de}. 

Tracing back results to their original source is a task on its own. We want to work on that issue more intensively in the 
future. If you observe possible enhancements in that direction, please let us know. Critique, suggestions for improvements
and feature requests are also highly welcome.

\section{Acknowledgement}
The authors want to thank the contributors Anna-Lena Horlemann-Trautmann, Ivan David Molina Naizir, Francesco Pavese, and Alexander Shishkin.

\iftoggle{arxiv}
{

}
{
  \bibliographystyle{plain}
  \bibliography{tables_subspace_codes}
}  
\pagebreak

\appendix
\footnotesize

\section{Tables for binary constant dimension codes}

\indent

\begin{tabular}{c|cc}
$n=6$ &  2 &  3 \\ 
\hline
4 &  21 *  $(131044)$  &  77 *  $(5)$ \\
6 & &  9 *  $(1)$ \\
\end{tabular}

\begin{tabular}{c|cc}
$n=7$ &  2 &  3 \\ 
\hline
4 &  41  &  333 - 381 \\
6 & &  17 *  $(715)$ \\
\end{tabular}

\begin{tabular}{c|ccc}
$n=8$ &  2 &  3 &  4 \\ 
\hline
4 &  85  &  1326 - 1493  &  4801 - 6477 \\
6 & &  34 $(\ge624)$  &  257 *  $(2)$ \\
8 & & &  17 *  $(8)$ \\
\end{tabular}

\begin{tabular}{c|ccc}
$n=9$ &  2 &  3 &  4 \\ 
\hline
4 &  169  &  5986 - 6205  &  37265 - 50861 \\
6 & &  73  &  1033 - 1156 \\
8 & & &  33 \\
\end{tabular}

\begin{tabular}{c|cccc}
$n=10$ &  2 &  3 &  4 &  5 \\ 
\hline
4 &  341  &  23870 - 24697  &  301213 - 423181 &  1178539 - 1678413 \\
6 & &  145  &  4173 - 4977  &  32923 - 38148 \\
8 & & &  65  &  1025 - 1089 \\
10 & & & &  33 \\
\end{tabular}

\begin{tabular}{c|cccc}
$n=11$ &  2 &  3 &  4 &  5 \\ 
\hline
4 &  681  &  97526 - 99717  &  2383041 - 3370315  &  18728043 - 27943597 \\
6 & &  290  &  16717 - 19785   &  263478 - 328641 \\
8 & & &  129 - 132  &  4097 - 4289 \\
10 & & & &  65 \\
\end{tabular}

\begin{tabular}{c|ccccc}
$n=12$ &  2 &  3 &  4 &  5 &  6 \\ 
\hline
4 &  1365  &  385515 - 398385  &  19673822 - 27222741 &  299769965 - 445207739  &  1212491081 - 1816333805 \\
6 & &  585  &  66862 - 79170  &  2105077 - 2613533  & 16865102 - 21361665 \\
8 & & &  273  &  16401 - 17436  &  262165 - 278785 \\
10 & & & &  129  &  4097 - 4225 \\
12 & & & & &  65 \\
\end{tabular}

\begin{tabular}{c|ccccc}
$n=13$ &  2 &  3 &  4 &  5 &  6 \\ 
\hline
4 &  2729  &  1597245 $(\ge512)$  &  157328094 - 217544769 &  4794061075 - 7192950693  &  38325127529 - 57884072859 \\
6 & &  1169  &  268130 - 319449  &  16835124 - 20918754  &  269057345 - 339800773 \\
8 & & &  545  &  65793 - 72131  &  2097225 - 2266956 \\
10 & & & &  257 - 259  &  16385 - 16769 \\
12 & & & & &  129 \\
\end{tabular}

\section{Tables for ternary constant dimension codes}
\indent

\begin{tabular}{c|cc}
$n=6$ &  2 &  3 \\ 
\hline
4 &  91  &  754 - 784 \\
6 & &  28 *  $(7)$ \\
\end{tabular}

\begin{tabular}{c|cc}
$n=7$ &  2 &  3 \\ 
\hline
4 &  271  &  6978 - 7651 \\
6 & &  82 \\
\end{tabular}

\begin{tabular}{c|ccc}
$n=8$ &  2 &  3 &  4 \\ 
\hline
4 &  820  &  60259 - 68374  &  543142 - 627382 \\
6 & &  244 - 248  &  6562 - 6724 \\
8 & & &  82 \\
\end{tabular}

\begin{tabular}{c|ccc}
$n=9$ &  2 &  3 &  4 \\ 
\hline
4 &  2458  &  550291 - 620740  & 14585908 - 16821712 \\
6 & &  757  &  59077 - 61010 \\
8 & & &  244 \\
\end{tabular}

\begin{tabular}{c|cccc}
$n=10$ &  2 &  3 &  4 &  5 \\ 
\hline
4 &  7381  &  5086975 - 5582305  &  394218370 - 458168194  &  3554720608 - 4104497728 \\
6 & &  2269  &   	532195 - 558739  &  14349660 - 14886440 \\
8 & & &  730 - 732  &  59050 - 59536 \\
10 & & & &  244 \\
\end{tabular}

\begin{tabular}{c|cccc}
$n=11$ &  2 &  3 &  4 &  5 \\ 
\hline
4 &  22141  &  45782788 - 50289022  &  10639658410 - 12361037515  &  286680643528 - 335382904522 \\
6 & &  6805 - 6809  &  4789947 - 5024299  &  387447165 - 409001563 \\
8 & & &  2188 - 2201  &  531442 - 535824 \\
10 & & & &  730 \\
\end{tabular}

\section{Tables for quaternary constant dimension codes}
\indent

\begin{tabular}{c|cc}
$n=6$ &  2 &  3 \\ 
\hline
4 &  273  &  4137 - 4225 \\
6 & &  65 \\
\end{tabular}

\begin{tabular}{c|cc}
$n=7$ &  2 &  3 \\ 
\hline
4 &  1089  &  66828 - 70993 \\
6 & &  257 \\
\end{tabular}

\begin{tabular}{c|ccc}
$n=8$ &  2 &  3 &  4 \\ 
\hline
4 &  4369  &  1054373 - 1132817  &  16874321 - 18245201 \\
6 & &  1025 - 1033  &  65537 - 66049 \\
8 & & &  257 \\
\end{tabular}

\begin{tabular}{c|ccc}
$n=9$ &  2 &  3 &  4 \\ 
\hline
4 &  17473  &  16947673 - 18179409 & 1078557605 - 1164549201 \\
6 & &  4161  &  1048641 - 1061929 \\
8 & & &  1025 \\
\end{tabular}

\begin{tabular}{c|cccc}
$n=10$ &  2 &  3 &  4 &  5 \\ 
\hline
4 &  69905  &  273727509 - 290821441  &  69038576145 - 74754799185  &  1105471620389 - 1193662931025 \\
6 & &  16641  &  16781353 - 17110273  &  1073745960 - 1088477225 \\
8 & & &  4097 - 4102  &  1048577 - 1050625 \\
10 & & & &  1025 \\
\end{tabular}

\begin{tabular}{c|cccc}
$n=11$ &  2 &  3 &  4 &  5 \\ 
\hline
4 &  279617  &  4379640165 - 4654011921 &  4418468947289 - 4783502911565  &  282679561437637 - 306494895880785 \\
6 & &  66561 - 66569  &  268502284 - 273715273 &  68719805936 - 70152169473 \\
8 & & &  16385 - 16418  &  16777217 - 16818202 \\
10 & & & &  4097 \\
\end{tabular}

\section{Table for (unrestricted) binary subspace codes}
\indent

\begin{tabular}{c|ccccccccc}
$q=2$ &  1 &  2 &  3 &  4 &  5 &  6 &  7 &  8 &  9 \\ 
\hline
1 &  2 *  $(1)$  &  5 *  $(1)$  &  16 *  $(1)$  &  67 *  $(1)$  &  374 *  $(1)$  &  2825 *  $(1)$  &  29212 *  $(1)$  &  417199 *  $(1)$  &  8283458 *  $(1)$ \\
2 & &  3 *  $(1)$  &  8 *  $(1)$  &  37 *  $(1)$  &  187 *  $(1)$  &  1521 *  $(1)$  &  14606 *  $(1)$  &  222379 *  $(1)$  &  4141729 *  $(1)$ \\
3 & & &  2 *  $(2)$  &  5 *  $(3)$  &  18 *  $(24298)$  &  108 - 117  &  614 - 776  &  5687 - 9191  &  71427 - 107419 \\
4 & & & &  5 *  $(1)$  &  9 *  $(7)$  &  77 *  $(4)$  &  334 - 388  &  4803 - 6479  &  37267 - 53710 \\
5 & & & & &  2 *  $(3)$  &  9 *  $(4)$  &  34 *  $(20)$  &  263 - 326  &  1996 - 2458 \\
6 & & & & & &  9 *  $(1)$  &  17 *  $(928)$  &  257 *  $(6)$  &  1034 - 1240 \\
7 & & & & & & &  2 *  $(4)$  &  17  &  65 - 66 \\
8 & & & & & & & &  17 *  $(7)$  &  33 \\
9 & & & & & & & & &  2 *  $(5)$ \\
\end{tabular}

\section{Table for (unrestricted) ternary subspace codes}
\indent

\begin{tabular}{c|ccccccccc}
$q=3$ &  1 &  2 &  3 &  4 &  5 &  6 &  7 &  8 &  9 \\ 
\hline
1 &  2 *  $(1)$  &  6 *  $(1)$  &  28 *  $(1)$  &  212 *  $(1)$  &  2664 *  $(1)$  &  56632 *  $(1)$  &  2052656 *  $(1)$  &  127902864 *  $(1)$  &  13721229088 *  $(1)$ \\
2 & &  4 *  $(1)$  &  14 *  $(1)$  &  132 *  $(1)$  &  1332 *  $(1)$  &  34608 *  $(1)$  &  1026328 *  $(1)$  &  77705744 *  $(1)$  &  6860614544 *  $(1)$ \\
3 & & &  2 *  $(2)$  &  10  &  56  &  764 - 966 &	13248 - 15394 &	544431 - 758228 &	29137055 - 34143770  \\
4 & & & &  10 *  $(2)$  &  28  &  754 - 784 &	6979 - 7696 &	543144 - 627384 &	14585910 - 17071886  \\
5 & & & & &  2 *  $(3)$  &  28  &  163 - 164  &  6574 - 7222  &  117621 - 123048 \\
6 & & & & & &  28 *  $(7)$  &  82  &  6562 - 6724 &	59078 - 61962  \\
7 & & & & & & &  2 *  $(4)$  &  82  &  487 - 488 \\
8 & & & & & & & &  82  &  244 \\
9 & & & & & & & & &  2 *  $(5)$ \\
\end{tabular}

\section{Table for (unrestricted) quaternary subspace codes}
\indent

\begin{tabular}{c|cccccccc}
$q=4$ &  1 &  2 &  3 &  4 &  5 &  6 &  7 &  8 \\ 
\hline
1 &  2 *  $(1)$  &  7 *  $(1)$  &  44 *  $(1)$  &  529 *  $(1)$  &  12278 *  $(1)$  &  565723 *  $(1)$  &  51409856 *  $(1)$  &  9371059621 *  $(1)$ \\
2 & &  5 *  $(1)$  &  22 *  $(1)$  &  359 *  $(1)$  &  6139 *  $(1)$  &  379535 *  $(1)$  &  25704928 *  $(1)$  &  6269331761 *  $(1)$ \\
3 & & &  2 *  $(2)$  &  17  &  130  &  4154 - 4771 &	131318 - 142313 &	16881731 - 20449159  \\
4 & & & &  17 *  $(3)$  &  65  &   	4137 - 4225 &	66829 - 71156 &	16874323 - 18245203  \\
5 & & & & &  2 *  $(3)$  &  65  &  513 - 514  &  65557 - 68117 \\
6 & & & & & &  65  &  257  &  65537 - 66049  \\
7 & & & & & & &  2 *  $(4)$  &  257 \\
8 & & & & & & & &  257 \\
\end{tabular}

\end{document}